\documentclass[a4paper]{amsart}
\usepackage[english]{babel}

\usepackage{latexsym}

\usepackage{enumerate}

\usepackage{graphicx}
\usepackage{graphics}
\usepackage{subfigure}
\usepackage{amsmath}
\usepackage{amsfonts}
\usepackage{amssymb}
\usepackage{amscd}
\usepackage{amsthm}
\usepackage{url}

\usepackage{multirow}
\usepackage{fancyvrb}
\usepackage{color}

\usepackage{fancyvrb}

\usepackage{listings}
\usepackage{xcolor}
\definecolor{mGreen}{rgb}{0,0.6,0}
\definecolor{mGray}{rgb}{0.5,0.5,0.5}
\definecolor{mPurple}{rgb}{0.58,0,0.82}
\definecolor{backgroundColour}{rgb}{0.95,0.95,0.92}

\lstdefinestyle{CStyle}{
	backgroundcolor=\color{backgroundColour},   
	commentstyle=\color{mGreen},
	keywordstyle=\color{magenta},
	numberstyle=\tiny\color{mGray},
	stringstyle=\color{mPurple},
	basicstyle=\footnotesize,
	breakatwhitespace=false,         
	breaklines=true,                 
	captionpos=b,                    
	keepspaces=true,                 
	numbers=left,                    
	numbersep=5pt,                  
	showspaces=false,                
	showstringspaces=false,
	showtabs=false,                  
	tabsize=2,
	language=C
}

\usepackage{hyperref}
\hypersetup{colorlinks=false}

\usepackage{pgfplots}
\tikzset{ font={\fontsize{9pt}{12}\selectfont}}

\newtheorem{theorem}{Theorem}[section]
\newtheorem{proposition}[theorem]{Proposition}
\newtheorem{lemma}[theorem]{Lemma}

\newtheorem{definition}[theorem]{Definition}
\newtheorem{remark}[theorem]{Remark}

\newcommand{\C}{\mathbb{C}}

\newcommand{\irrat}{\mathbb{R}\setminus\mathbb{Q}}

\newcommand{\wcom}{\widehat{\mathbb{C}}}

\newcommand{\rat}{\mathbb{Q}}

\begin{document}

\title{Computing  parameter planes of iterative root-finding methods with several  free critical points }
\begin{author}[B.~Campos]{Beatriz Campos}
	\email{campos@uji.es}
	\address{ %
		Instituto de Matem\'aticas y Aplicaciones de Castell\'on, Universitat Jaume I, Spain. ORCID: \url{https://orcid.org/0000-0001-9205-0256}}
\end{author}

\begin{author}[J.~Canela]{Jordi Canela}
	\email{canela@uji.es}
	\address{ %
		Instituto de Matem\'aticas y Aplicaciones de Castell\'on, 
		Universitat Jaume I,  Spain. ORCID: \url{https://orcid.org/0000-0001-7879-5438}}
\end{author}

\begin{author}[A.~Rodríguez-Arenas]{Alberto Rodríguez-Arenas}
	\email{arenasa@uji.es}
	\address{Universitat Poilitècnica de València, Spain. ORCID: \url{https://orcid.org/0000-0002-4571-2031}}
\end{author}

\begin{author}[P.~Vindel]{Pura Vindel}
	\email{vindel@uji.es}
	\address{ %
		Instituto de Matem\'aticas y Aplicaciones de Castell\'on,
		Universitat Jaume I,  Spain.  ORCID: \url{https://orcid.org/0000-0001-8422-4738}}
\end{author}

\maketitle

\begin{abstract}
In this paper we present an algorithm to obtain the parameter planes of families of root-finding methods with several free critical points.  The parameter planes show the joint behaviour of all critical points.
This algorithm avoids the inconsistencies arising from the relationship between the different critical points  as well as
the indeterminacy  caused by the square roots involved in their computation.

We analyse the suitability of this algorithm by drawing the parameter planes of different Newton-like methods with two and three critical points.
We also present some results of the expressions of the Newton-like operators and their derivatives  in terms of palindromic polynomials, and we show how to obtain the expression of the critical points of a Newton-like method  with real coefficients.



\end{abstract}

\section{Introduction}

Iterative root-finding methods are used to solve equations whose solutions cannot be obtained by means of algebraic procedures. The development of new root-finding methods has become a very active area of research: it is sought to find new methods which increase the order of convergence to the solutions of the equation and have better computational efficiency. However, the radii of convergence may decrease as the order of the methods increases. At this point, a dynamical study can provide valuable information on the behaviour  of these methods in a qualitative way. 


When we apply an iterative root-finding method to solve the non-linear equation $f(z)=0$ we obtain an operator $O_f$. If $z^*$ is a solution of the equation and an initial guess $z_0$ is close enough to $z^*$, the operator provides a sequence 
$$z_{k+1}=O_f(z_k), \;\; k\geq0,$$
that converges to $z^*$. The sequence $\{z_k\}_{k\geq 0}$ defines the orbit of the point $z_0$. Nonetheless, in general we cannot know a priori if an initial guess $z_0$ is close  enough to a root. A dynamical study can provide information on the kind of asymptotic behaviour presented by the orbits of an initial condition $z_0$. In particular, it can detect whether there are attracting cycles or other stable behaviour not corresponding to the basins of attraction of the solutions of $f(z)=0$. Such stable behaviours would provide open sets of initial conditions which do not converge to any of the roots, which is an important drawback when applying the algorithm. In many cases, if $f$ is a rational map, then $O_f$ is also a rational map and we can use the theory of complex dynamics (see \S~\ref{subsec:introholo}) to study the existence of such stable domains. Indeed, all stable behaviours are related to a critical point, a point $c$ such that $O'_f(c)=0$, so it is enough to study the orbits of all critical points of $f$.

If we study a family of root-finding algorithms depending on parameters (or if $f$ depends on parameters), the operator $O_f$ will also depend on parameters. In this scenario, parameter planes play an important role on helping us understand the family.  For each parameter, we can use the orbit of critical points to determine if there are stable domains other than the basins of attraction of the roots and then plot the parameter accordingly.   Using the parameter plane we can find the members of the family with better behaviour.

When the operator has a single free critical point (a critical point which is not fixed under the dynamics of $O_f$), each colour of the parameter plane explains the asymptotic behaviour of that critical point (see \cite{CCV-gradon},  \cite{CTV-gato}, for example). However, as the order of convergence of the algorithm increases, the number of free critical points also increases, and drawing parameter planes becomes challenging (see, for instance, \cite{operador6+4},  \cite{argyros}, \cite{AM-pg5},  \cite{tresCrit}, \cite{Behl}, \cite{CCTV-varianteCheby}, \cite{CHUN2012}, \cite{CTV-ermakov}, \cite{Geum},     \cite{sharma},  \cite{ZCT-pg8} and references therein). A usual approach to tackle this problem is to produce a different parameter plane for each different free critical point. When considered simultaneously, these parameter planes provide complete information of the asymptotic behaviour of all critical orbits. However, when regarded separately these planes may contain inconsistencies which are usually due to changes of determination of roots which appear in the definition of the critical points and the dynamical relations amongst them (see Figure~\ref{c1c3OP2v2}, upper). Moreover, some bifurcations may be difficult to understand when plotting parameter planes separately. For instance, if one critical point is captured by an attracting cycle controlled by another critical point (a capture parameter), the parameter planes may seem unusual when regarded separately (see, for instance, Figure~\ref{fig:param3}). An alternative approach to this problem is to plot a single parameter plane which considers the dynamics of these critical points simultaneously, plotting the parameter in black if any of the critical orbits does not converge to the roots. Those black parameters correspond to operators for which the root-finding algorithm may not behave appropriately: there can be stable behaviours other than the roots. This approach is followed, for instance, in \cite{Gutierrezcubicos} and \cite{CEGJ}. However, those parameter planes may have the disadvantage of losing the information of how many critical orbits fail to converge to the roots. For instance, in \cite[Figure 5]{CEGJ} there are black parameters for which only one critical orbit fails to converge to the roots and other parameters for which no free critical orbit converges to the roots. This information is relevant since we can have as many attracting cycles not coming from the basins of attraction of the roots as critical orbits failing to converge to the roots.

The goal of this paper is to present an algorithm to draw parameter planes taking into account all critical points simultaneously in a single plane and not losing any information. This algorithm is presented in \S~\ref{algoritmo} and is based on the escaping algorithm and can be used with no modification for any number of free critical points. The idea of the algorithm is the following. If all critical points converge to the roots, then we use a scaling of colours which indicates the slowest time of convergence to a root amongst all critical points. This criterion avoids analysing to which roots the critical orbits converge, since changes on the determination of the roots could lead to lines in the parameter planes which do not actually correspond to bifurcations, similar to what can be observed in the upper planes of Figure~\ref{c1c3OP2v2}. If any of the critical orbits does not converge to a root, then we plot the parameter with a different colour depending on the number of critical points which do not escape (see Figures~\ref{c1c3OP2v2}, \ref{fig:param3} and \ref{fig:param3crit}). Along the paper we explain how using this algorithm can help to better interpret the bifurcations in parameter plane. We also present different modifications of the algorithm that can help us get extra information (see Figures~\ref{fig:param3-2col} and \ref{fig:paramO4}). The implementation in C of the algorithm is available upon request.

Even though the algorithm presented can be adapted to plot the parameter plane of any family of root-finding algorithms, the implementation that we present is done keeping in mind the so called \textit{Newton-like methods}. These methods are variations of Newton's method and are analysed in \cite{Newton-like}. Many of the root-finding algorithms in the literature are Newton-like methods  (\cite{AM-pg5},  \cite{CCV-gradon+k},  \cite{CCTV-familiac},  \cite{CV-OstrowskiChun},    \cite{CFMT-pg3}, \cite{CGTVV-king},  \cite{CGMT-pg4}, \cite{CMQT-pg1}, \cite{CTV-gato}, \cite{CTV-ermakov},   \cite{MA-pg2},   \cite{ZCT-pg8}, for example). When  Newton-like algorithms are applied on quadratic polynomials $p(z)=z^2-c$, an intrinsic symmetry appears in the operator obtained. We prove that, after applying a conjugacy that sends the roots to $z=0$ and $z=\infty$, such operators have the following generic expression:

\begin{equation*}
    O(z)=z^n \frac{a_k+a_{k-1} z+...+a_{1} z^{k-1}+ z^k}{1+a_{1} z+...+a_{k-1} z^{k-1}+a_k z^k}=z^n \prod_{i=1}^{k}
\frac{(z-r_i)}{(1-r_i z)}.
\end{equation*}

Actually, in \cite{Newton-like} we study such methods when applied on polynomials  $p(z)=z^d-c$ and  we show that the maps obtained  are symmetric with respect to a rotation by a $d$th root of the unit.  The previous operator is obtained when restricting to $d=2$ and applying the  conjugacy, regardless of the method used. After applying the conjugacy, the operator $O(z)$ obtained is symmetric with respect to the map $z\rightarrow 1/z$. This symmetry is taken into consideration when implementing the algorithm. Indeed, in order to avoid inconsistencies in the colour scheme used when all critical orbits converge to the roots, we need to implement stop conditions for the convergence to $z=0$ and $z=\infty$ which are preserved by $z\rightarrow 1/z$. Moreover, if $c$ is a critical point of $O(z)$, by symmetry, then $\tilde{c}=1/c$ is also a critical point of $O(z)$ and their orbits have symmetric asymptotic behaviour. Therefore, we count each pair $\{c,1/c\}$ as a single free critical point and only iterate one of them when drawing parameter planes.

Up to this moment we have not talked about another crucial procedure to draw parameter planes. We need to actually compute the expressions of all critical points of the operator $O(z)$. As the degree of $O(z)$ increases, the number of critical points also increases (a rational map of degree $d$ has $2d-2$ critical points counting multiplicity), so obtaining expressions of all the critical points can be challenging. However, the operators $O(z)$ coming from Newton-like methods satisfy certain properties which may help us find all of their critical points. In  \S~\ref{estudioCriticos} we prove that the derivative of these operators gives rise to palindromic polynomials, which allows us to halve the degree of the polynomial we need to solve. 
In particular, we stablish in Proposition \ref{proposicion} that 
the free critical points from the operator $O(z)$  with real coefficients satisfy that they are either  pairs of inverse real roots, or complex conjugates lying on the unit circle or  a set of four related roots of a quartic palindromic polynomial. The methods introduced in \S~\ref{estudioCriticos} are later used in \S~\ref{ejemplos} in order to obtain all critical points of different Newton-like root finding algorithm's for which we later  plot the parameter planes.

This paper is organized as follows. In \S~\ref{subsec:introholo} we finish the introduction by recalling the basic concepts of complex dynamics. In  \S~\ref{estudioCriticos} we present the relation of the operators $O(z)$ obtained applying Newton-like root finding algorithm's to quadratic polynomials and prove that the numerator of their derivative is palindromic. Afterwards we show different techniques to solve palindromic polynomials in order to be able to obtain the critical points of those operators.
  In \S~\ref{algoritmo} we explain in detail the algorithm for drawing parameter planes. Next, in \S~\ref{ejemplos} we illustrate the convenience of using the algorithm by drawing the parameter planes of different Newton-like methods with two and three critical points. We also plot several dynamical planes corresponding to each of the examples for a better understanding of the colours in the parameter plane. Moreover, we explain different modifications to the algorithm that can be implemented to obtain a better understanding of the parameter planes. 
\subsection{Introduction to complex dynamics}\label{subsec:introholo}

For a better understanding of the exhibited  results, we recall some basic concepts of complex dynamics. For a more detailed introduction to the topic we refer to \cite{Beardon, Milnor}.

Given a rational map $Q:\wcom\rightarrow\wcom$, where $\wcom$ denotes the Riemann sphere, we consider the dynamical system provided by the iterates of $Q$. A point $z_0$ is called fixed  if $Q(z_0)=z_0$. A point $z_0$ is called periodic of period $p\geq 1$ if $Q^p(z_0)=z_0$ and $Q^\ell(z_0)\neq z_0$ for all $\ell<p$. In the later case we denote by $\langle z_0\rangle =\{z_0, z_1,\ldots,z_{p-1}\}$, where $z_{\ell}=Q(z_{\ell-1})$, the cycle of period $p$ generated by $z_0$. The multiplier of a fixed point is given by $\lambda(z_0)=Q'(z_0)$.
 Similarly, the multiplier of a periodic point is given by $\lambda\left(\langle z_0\rangle\right)=\left(Q^p\right)'(z_0)=Q'(z_0)\cdot \ldots \cdot Q'(z_{p-1})$. A periodic or fixed point $z_0$ is called \textit{attracting} (resp.\ \textit{superattracting}) if $|\lambda|<1$ (resp.\ $\lambda=0$), \textit{repelling} if $|\lambda|>1$, and \textit{indifferent} if $|\lambda|=1$. An indifferent point $z_0$ is called \textit{parabolic} (or \textit{rationally indifferent}) if $\lambda=e^{2\pi i r/s}$ with $r/s\in\rat$. If $\lambda=e^{2\pi i \theta}$ with $\theta\in\irrat$ the point $z_0$ is called \textit{irrationally indifferent}. 
 Attracting fixed (or periodic) points $z_0$ have associated a basin of attraction $\mathcal{A}(z_0)$ associated to them, which consists of the set of points that converge to $z_0$ (or the cycle $\langle z_0\rangle$) under iteration of $Q$. 
Similarly, the basin of attraction $\mathcal{A}(z_0)$ of an a parabolic fixed (or periodic) point is defined as the set of points which converge to $z_0$ (or $\langle z_0\rangle$). Unlike in the attracting case, a parabolic point $z_0$ belongs to the boundary of $\mathcal{A}(z_0)$, $z_0\notin \mathcal{A}(z_0)$. With respect to the irrationally indifferent point $z_0$, if the map $Q$ ($Q^p$ in the periodic case) is conjugate to the rigid rotation $z\rightarrow \theta\cdot z$ in some neighbourhood of $z_0$ we say that $z_0$ is a \textit{Siegel point} and the maximal domain of the conjugation is called \textit{Siegel disk}. Otherwise we say that $z_0$ is a \textit{Cremer point}.

The iteration of $Q$ defines a completely invariant partition of $\wcom$. The \textit{Fatou set} $\mathcal{F}(Q)$ is defined as the set of points $z\in\wcom$ such that the family of iterates of $Q$ is normal in some open neighbourhood of $z$. Its complement, the \textit{Julia set} $\mathcal{J}(Q)=\wcom\setminus\{\mathcal{F}(Q)\}$, is closed and corresponds to the set of points with chaotic behaviour. The connected components of $\mathcal{F}(Q)$ are called Fatou components and are mapped under iteration of $Q$ amongst themselves. It follows from Sullivan's No Wandering Theorem \cite{Su} that all Fatou components of a rational map are either periodic or preperiodic. All periodic Fatou components of a rational map are either basins of attraction of attracting or parabolic cycles, or simply connected rotation domains (Siegel disks) or doubly connected rotation domains (Herman rings). Moreover, all these periodic Fatou components are related to a critical point, i.e.\ a point $c\in \wcom$ such that $Q'(c)=0$. Indeed, all attracting and parabolic basins of attraction contain, at least, a critical point. Furthermore, given any Siegel disk or Herman ring $U$ there is a critical point (two in the case of Herman rings) whose orbit accumulates on $\partial U$. If a critical point is not a fixed point of $Q$ it is called free critical point. Two o more critical points can satisfy relations among them that imply a symmetry in their dynamics; 
therefore, in order to detect all stable behaviours of the map $Q$ it is enough to study the asymptotic behaviour of all free critical points of $Q$ up to symmetry.

\section{Critical points of Newton-like methods} \label{estudioCriticos}
In this section, we study properties of the operators obtained when  applying Newton-like methods to polynomials of degree two that allow us to obtain the expressions of all critical points in terms of the parameter. As proved in \cite{Newton-like}, these operators have the generic expression:

\begin{equation}\label{op}
    O(z)=z^n \frac{a_k+a_{k-1} z+...+a_{1} z^{k-1}+ z^k}{1+a_{1} z+...+a_{k-1} z^{k-1}+a_kz^k}
\end{equation}
with $a_k \neq 0$.

We can observe that  the polynomials in the numerator and denominator of this expression have the same coefficients in reciprocal order.

\begin{definition} Two degree $n$ polynomials $p(z)=a_0+a_1 z+...+a_{n-1}z^{n-1}+ a_n z^n$ and $q(z)=b_0+b_1 z+...+b_{n-1}z^{n-1}+ b_n z^n$ are  reciprocal if  $b_i= a_{n-i}$, for $0 \leq i \leq n$.
    \end{definition}

In the following, we study some results involving reciprocal polynomials.

  \begin{lemma}\label{simetria}
The quotient of two reciprocal polynomials p and $\widehat{p}$ satisfies the symmetry property
\[
\frac{p(z)}{\widehat{p}(z)}=\frac{1}{\frac{p(1/z)}{\widehat{p}(1/z)}}.
\]
\end{lemma}

\begin{proof}
      Let us consider a degree $n$ polynomial $p(z)=a_0+a_1 z+...+a_{n-1}z^{n-1}+ a_n z^n$. The reciprocal polynomial of $p(z)$ is $\widehat{p}(z)=a_n+a_{n-1} z+...+a_{1}z^{n-1}+ a_0 z^n$, that can be written as $\widehat{p}(z)=z^n p(1/z)$. Then,
\[    
\frac{p(1/z)}{\widehat{p}(1/z)}  =\frac{z^{-n}\widehat{p}(z)}{(1/z)^n p(z)}=\frac{\widehat{p}(z)}{ p(z)}.
\]    
   \end{proof}
    
\begin{remark}    
This property implies  that the operator (\ref{op}) satisfies $O(z)=\frac{1}{O(1/z)}$. So, the strange fixed points different from $z=1$ and $z=-1$  and the critical points different from $z=0$ of the operator $O(z)$ come in inverse pairs.
\end{remark}

The fact that reciprocal polynomials appear in the expression of the operator $O(z)$ leads to a special type of polynomials in the expression of $O'(z)$, the so-called palindromic polynomials.

 \begin{definition} A polynomial $P(z)=A_0+A_{1} z+...+A_{n-1}z^{n-1}+ A_n z^n$ is called palindromic  if $A_{i}=A_{n-i},$ $0 \leq i \leq n$.
    \end{definition}

   \begin{definition} A polynomial $P(z)=A_0+A_{1} z+...+A_{n-1}z^{n-1}+ A_n z^n$ is called antipalindromic  if $A_{i}=-A_{n-i},$ $0 \leq i \leq n$.
    \end{definition}

Some properties of this type of polynomials are given in the following proposition (see, for example, \cite{palindromo} for a proof of the result):
\begin{proposition}\label{prop:proerties palindromic}
	Palindromic and antipalindromic polynomials satisfy the following properties:
\begin{enumerate}[a)]
  \item  If $\alpha$ is a root of a polynomial that is either palindromic or antipalindromic, then $1/\alpha$ is also a root and
has the same multiplicity.
  \item The converse is true: if a polynomial is such that if $\alpha$ is a root then $1/\alpha$ is also a root of the same multiplicity, then the polynomial is either palindromic or antipalindromic.
  \item The sum of two palindromic (antipalindromic) polynomials is a palindromic (antipalindromic) polynomial.
  \item The product of a constant by a palindromic (antipalindromic) polynomial is a palindromic (antipalindromic) polynomial.
  \item The product of two palindromic  or two antipalindromic polynomials is palindromic.
  \item A palindromic polynomial $P(z)$ of odd degree is a multiple of $z+1$ (it has -1 as a root) and its quotient by $z+1$ is also palindromic.
  \item An antipalindromic polynomial $Q(z)$ is a multiple of $z-1$ (it has 1 as a root) and its quotient by $z-1$ is palindromic.
\end{enumerate}
\end{proposition}

Next we prove that the critical points of a Newton-like method applied on degree $2$ polynomials are the roots of a palindromic polynomial.

 \begin{lemma}\label{product}
The product of two reciprocal polynomials is a palindromic polynomial.
\end{lemma}

    \begin{proof}
    Let us consider a polynomial $p(z)=a_0+a_1 z+...+a_{n-1}z^{n-1}+ a_n z^n$ of degree n. The reciprocal polynomial of $p(z)$ is $\widehat{p}(z)=a_n+a_{n-1} z+...+a_{1}z^{n-1}+ a_0 z^n$, that can be written as $\widehat{p}(z)=z^n p(1/z)$. Then,

    \[
p(z)\widehat{p}(z)=\sum_{i=0}^{n} a_i z^i \cdot z^n \sum_{j=0}^{n} a_j (\frac{1}{z})^j=\sum_{i,j=0}^{n} a_i a_j z^{n+i-j}=\sum_{k=0}^{2n} b_k z^k
\]

where:

\begin{eqnarray*}
  b_k &=& \sum_{i=0}^{k} a_i a_{n-k+i}, \ \ 0\leq k \leq n\\
  b_k &=& \sum_{i=0}^{2n-k} a_{k-n+i} a_{i}, \ \ n\leq k \leq 2n.
\end{eqnarray*}

 Let us prove that  $b_{k}=b_{2n-k}$, for  $0\leq k \leq n$:

 \begin{eqnarray*}
  b_{2n-k}=\sum_{i=0}^{2n-(2n-k)} a_{(2n-k)-n+i} a_{i}=\sum_{i=0}^{k} a_{n-k+i} a_{i}=b_{k}.
\end{eqnarray*}

We can conclude that the product $p(z) \widehat{p}(z) $ is a palindromic polynomial.
\end{proof}

In \cite{Newton-like} we deduce that $z=1$ is a fixed point of  the operator $O(z)$. Moreover, if $n+k$ is odd, then $z=-1$ is also a fixed point. As we have remarked above, from Lemma \ref{simetria} it is easy to prove that the strange fixed points of the operator $O(z)$, different from $z=1$ and $z=-1$, come in inverse pairs. The same occurs for the critical points different from $z=0$. In fact, let us see that these critical points are the roots of a palindromic polynomial.

 \begin{lemma}
The polynomial in the numerator of $O^{\prime}(z)$ is palindromic.
\end{lemma}

\begin{proof}

Let us consider the operator $O(z)$ given in (\ref{op}) written as 
\[
 O(z)=z^{n} \frac{p(z)}{ \widehat{p}(z)}.
\]

Then, we have that:

\[
 O^{\prime}(z)=z^{n-1} \frac{n \widehat{p}(z)p(z)+  z \left( \widehat{p}^{\prime}(z)p(z)-\widehat{p}(z) p^{\prime}(z)  \right) }{p^2(z) }
\]
and the critical points are $z=0, \infty$ and the roots of the polynomial $P(z)=n \widehat{p}(z)p(z)+  z \left( \widehat{p}^{\prime}(z)p(z)-\widehat{p}(z) p^{\prime}(z)  \right) $ that appears in the numerator of $O^{\prime}(z)$. Let us see that $P(z)$ is a palindromic polynomial.

From the expressions $p(z)$ and $ \widehat{p}(z) $ we obtain $p'(z)$ and $ \widehat{p'}(z) $:
\begin{eqnarray*}
 p(z)=\sum_{i=0}^{n} a_i z^i &\Rightarrow&   p^{\prime}(z)=\sum_{i=0}^{n}i a_i z^{i-1} \\
\widehat{p}(z)=\sum_{j=0}^{n} a_j z^{n-j} &\Rightarrow&\widehat{p}^{\prime}(z)=\sum_{j=0}^{n}(n-j) a_j z^{n-j-1},
\end{eqnarray*}
and we can write:
\begin{eqnarray*}
  z \left(p(z)\widehat{p}^{\prime}(z)-p^{\prime}(z)\widehat{p}(z) \right) &=& z \left( \sum_{i=0}^{n}i a_i z^{i-1} \sum_{j=0}^{n} a_j z^{n-j}       \right)- z \left( \sum_{i=0}^{n} a_i z^i \sum_{j=0}^{n}(n-j) a_j z^{n-j-1}      \right) \\
    &=& \sum_{i,j=0}^{n}(i-n+j)a_i a_j z^{i+n-j}=  \sum_{k=0}^{2n}b_k z^{k},
\end{eqnarray*}
where:

\begin{eqnarray*}
  b_k &=& \sum_{i=0}^{k}(k-2i) a_i a_{n+i-k}, \ \ 0\leq k \leq n\\
  b_k &=& \sum_{i=0}^{2n-k}(2n-2i-k) a_{i+k-n} a_{i}, \ \ n\leq k \leq 2n.
\end{eqnarray*}

We can check that the coefficients verify $b_k=b_{2n-k}$, for $0\leq k \leq n$:
\begin{eqnarray*}
b_{2n-k} &=& \sum_{i=0}^{2n-(2n-k)}(2n-2i-(2n-k))  a_{i+(2n-k)-n} a_i=\sum_{i=0}^{k}(k-2i)  a_{n+i-k} a_i=  b_k.
\end{eqnarray*}

Then, the expression $z \left(p(z)\widehat{p}^{\prime}(z)-p^{\prime}(z)\widehat{p}(z) \right) $ is a palindromic polynomial. By Lemma \ref{product}, we have that $\widehat{p}(z)p(z)$ is also palindromic. By applying properties of palindromic polynomials (Proposition~\ref{prop:proerties palindromic}),  we conclude that $P(z)$ is a palindromic polynomial. 
\end{proof}

Then, the critical points different from $z=0$ and $z=\infty$ are the roots of a palindromic polynomial. Moreover, if the coefficients of the palindromic polynomial in the numerator of $O^{\prime}(z)$ are real, it can be decomposed as a finite product of polynomials of degree at most four. To achieve this goal, we rely on the following theorem concerning palindromic polynomials (see \cite{barbeau}, for example):

\begin{theorem}
For a polynomial  $P(z)=a_0+a_1 z+...+a_{n-1}z^{n-1}+ a_n z^n$ with coefficients in $\C$ and degree $n$, the following conditions are equivalent:
\begin{itemize}
\item the polynomial has palindromic coefficients: $a_k=a_{n-k}$ for all $k$,
\item  $z^n P(1/z) = P(z)$,
\item  (if $n=2m$) $P(z) = z^m q(z + 1/z)$ for a polynomial $q$ with coefficients in $\C$ and degree $m$.
\end{itemize}
\end{theorem}

 Given that an odd degree  palindromic polynomial can be written  as $(z+1)$ multiplied by a palindromic polynomial of even degree, we restrict our study to
polynomials of even degree. Moreover, let us notice that the change of variable $x=z+\frac{1}{z}$  transforms a palindromic polynomial $p(z)$ of degree $2m$ into a polynomial $q(x)$ of degree $ m$.

In the following result we show how to find all the critical points of the operator (\ref{op}) when  the coefficients of the rational function are real.

\begin{proposition} \label{proposicion}
The free critical points of the rational function given in (\ref{op}) with real coefficients satisfy that they are either  pairs of inverse real roots, or complex conjugates and lie on the unit circle or they are a set of four related roots of a quartic palindromic polynomial.
\end{proposition}  

\begin{proof}  
By applying the fundamental theorem of algebra,  a polynomial $q(x)$ with real coefficients can be decomposed as a product of monomials (corresponding to their real roots) and quadratic polynomials (corresponding to their complex conjugate roots). 
 
From the above results, it is obtained that every palindromic polynomial  $P(z)$  with real coefficients can be factorized into a product of palindromic polynomials of order two and four:
\[
P(z)= K \prod_{i} (z^2+a_i z+1)     \prod_{j} (z^4+b_j z^3+ c_j z^2+ b_j z+1).
\]

This statement is easy to see since, from the previous theorem, if $P$ has $n=2m$ degree, it can be written as $P(z) = z^m q(z + 1/z)$ and,  by applying the fundamental theorem of algebra on $ q(z + 1/z)$ we obtain:

\begin{eqnarray*} 
q(z + 1/z )  &=&  K \prod_{i=1}^{m_1} \left( (z+1/z)+A_i \right) \prod_{j=1}^{m_2} \left ((z+1/z)^2+B_j (z+1/z)+ C_j\right ) \\
&=&  K \prod_{i=1}^{m_1} \frac{1}{z}\left (z^2+A_i z+1\right ) \prod_{j=1}^{m_2} \frac{1}{z^2} \left (z^4+B_j z^3+(2+C_j)z^2 +B_j z+1\right )
\end{eqnarray*}  
where $m=m_1+2m_2$.

When studying the solutions of the polynomials $P$ and $q$ we want to highlight the following considerations:

\begin{itemize}
\item Polynomials $(z^2+A_i z+1)$ can be decomposed as a product of two monomials when $|A_i| \geq 2$; so, the corresponding roots of $P(z)$ are real and inverse. If $|A_i| <2$,  the corresponding roots of $P(z)$ are complex conjugate and they are on the unit circle.

\item The roots  of the polynomial $\left ((z+1/z)^2+B_j (z+1/z)+ C_j\right )$ for
 $(z+1/z)$ are real for  $ B_j^2-4C_j \geq 0$; so, the corresponding roots of $P(z)$ are  complex conjugate and they are on the unit circle. When $B_j^2-4C_j<0$, the inverse of a root of $P(z)$  is not its conjugate, then it must be one
of a set of four related roots that satisfy a quartic palindromic polynomial. 

\item Moreover, it is easy to check that the roots of $P$ on the unit circle, considered as pairs of reciprocals, correspond to the roots of $q$ in the interval $[-2,2]$.
\end{itemize}
\end{proof}

So, from the above results, it follows that if we are able to obtain the corresponding decomposition, we can always find all the critical points of a rational function of the type given in (\ref{op}).

\section{The algorithm } \label{algoritmo}

When a family of rational maps has more than one free critical point, understanding the parameter plane can be tricky. In this case, a usual procedure is to plot the parameter plane of every critical point separately. However, this poses two problems. First, it might be challenging to understand the whole bifurcation locus by observing the different plots separately (see Figure~\ref{c1c3OP2v2}, upper plots). Indeed, many times the changes in the asymptotic behaviour of one critical orbit may be determined by changes on other orbits (for instance, when one critical orbit is captured by an attracting cycle ``controlled'' by another critical orbit), but it is difficult to understand such behaviour by observing the plots obtained by iterating each free critical point separately. Moreover, changes in the determination of the roots may lead to non-continuous parametrizations of the critical points (when the determination changes the different critical points ``permute'' amongst one other). These phenomenon leads to curves in the plots which may be confused with bifurcations (see Figure~\ref{c1c3OP2v2} and \ref{fig:paramO4}).

In order to avoid these problems we plot the parameter plane by studying the orbits of all free critical points simultaneously. The algorithm used works as follows. First we create a grid of points. Each point of the grid is associated to a parameter in the region of the parameter plane that we want to draw.  Then, for each of these parameters we compute all different free critical points. 

Once all critical points are defined, each of them is iterated up to a  given maximum number of iterates.  Upon each iterate we verify if the orbit has converged to any of the roots. Since for this paper we consider root-finding algorithms applied to quadratic polynomials and we conjugate the operator obtained so that the roots are placed at $0$ and $\infty$, in order to verify if we have convergence, we check if the iterate $z$ satisfies $|z|<eps$ (convergence to 0) or if $|z|>esc$. We use $esc=10^4$ and $eps=1/esc$.
 If the critical orbit converges to one of the roots before that, we stop the process and iterate the next critical point. Moreover, we store the information of the amount of iterates needed by the ``slowest'' critical point to converge to the roots.
 
If at the end of the process all free critical orbits converge to the roots, we plot the pixel using a scaling of colours which goes from red, to yellow, to pallid green, to blue, and up to white. If any of the critical orbits does not converge to the roots then we plot the pixel with a different colour depending on the number of orbits which do not converge to the roots (black if no free critical points converges to the roots, pink if only one critical point converges to the roots, dark green if two converge to the roots, etc.).

Even though in this paper we only apply the algorithm to families with up to three free critical orbits, the algorithm is designed to handle without modification any number of free critical orbits.

We want to make a remark on how the parameters $esc$ and $eps$ that are used to determine convergence to the roots are chosen. In this paper we work with Newton-like families applied to quadratic polynomials and, hence, the operators obtained are symmetric with respect to the map $z\rightarrow 1/z$ (see \cite{Newton-like}). It follows that if $c$ is a critical point then $1/c$ is also a critical point and their orbits are symmetric. In the program we only iterate one of each pair of critical points since both provide the same information. However, in order to guarantee that the information of how fast a critical point converges to the roots does not depend on the critical point chosen, the stop criterium needs to respect the symmetry $z\rightarrow 1/z$. This is why we choose $eps=1/esc$. 

\section{ Newton-like methods with more than one free critical point} \label{ejemplos}

In this section we apply the program to plot parameter planes of different operators with more than one free critical orbit obtained from applying Newton-like methods  on degree two polynomials.  Examples with more than one free critical point appear with some frequency in the literature, especially when high-order numerical methods are studied. We consider some cases that are representative of the type of dynamics that they give rise to.

As proved in \cite{Newton-like}, $z=1$ is always a strange fixed point of this type of methods; $z=-1$ is also a fixed point when $n+k$ is odd and it is a preimage of $z=1$ when $n+k$ is even. We check these statements in the methods that we study. 

As the operators have two or more critical points, if we plot separately the parameter plane of each of them there appear inconsistencies produced by the indeterminacies of the square roots in the expression of the critical points, as we can see in the figures that appear in each subsection. However, this problem does not appear when using the program since the parameter plane that it plots takes into account simultaneously the behaviour of all free critical orbits. 

The first three  subsections correspond to Newton-like systems with operators satisfying that $n \geq k$; we can observe that their parameter planes are similar to the parameter planes with one free critical point obtained in \cite{Newton-like}. 

The family studied in \S~\ref{seccion_Ermakov-Kalitkin family}  corresponds to  a limit case of this type of numerical methods, so it deserves a more detailed study. 
Although Ermakov-Kalitkin family has two free critical points, one of them is necessarily in the basin of attraction of the point $z=-1$, which is a parabolic point located on the boundary of two attractor petals. So, it can be considered as a family with a single free critical point.

Finally, in \S~\ref{seccion_tres criticos} we consider an example where we use our algorithm to obtain the parameter plane of a family with three free critical points.

\subsection{Fourth-order methods derived from the Kim family}\label{subsec:op2}

In \cite{CCT2013}, the authors study a parametric family of fourth-order methods coming from the Kim family. After applying it on quadratic polynomials they obtain the following operator:
\begin{equation}\label{op2}
    O_a(z)=\frac{z^4(1-a+4z+6z^2+4z^3+z^4)}{1+4z+6z^2+4z^3+(1-a)z^4}
\end{equation}
whose derivative is:
\begin{equation}\label{dop2}
    O^{\prime}_a(z)=\frac{4z^3(1+z)^4(-1+a+(-4-a)z+(-6+a)z^2+(-4-a)z^3+(-1+a)z^4)}{(1+4z+6z^2+4z^3+(1-a)z^4)^2}.
\end{equation}

The fixed points of $O_a(z)$ are $z=0$, $z=\infty$ and the solutions of the equation:
\[
(z-1)(1 + 5 z + 11 z^2 +  (14 + a)  z^3 + 11  z ^4 + 5  z ^5 + z ^6)=0.
\]
So,  there exist seven strange fixed points, $z=1$ and the six roots of the polynomial of degree six above. 
As proved in \cite{Newton-like}, the fixed point $z=1$ is attractive outside the circle
\[
|a-16|=64
\]
and the point $z=-1$ is a preimage of $z=1$.

The free critical points are the roots of the four-degree polynomial in the numerator of (\ref{dop2}). With the change $x=z+\frac{1}{z}$, the problem is reduced to find the solutions of the equation
\[
(1-a)x^2+(4+a)x+4+a=0.
\]
By undoing this change, the four free critical points are:
\begin{eqnarray*}
  c_1(a) &=& \frac{4+a-\sqrt{5a(4+a)}-\sqrt{10a(6-a)-2(4+a)\sqrt{5a(4+a)}}}{4(a-1)}, \\
  c_2(a) &=& \frac{4+a-\sqrt{5a(4+a)}+\sqrt{10a(6-a)-2(4+a)\sqrt{5a(4+a)}}}{4(a-1)}, \\
  c_3(a) &=& \frac{4+a+\sqrt{5a(4+a)}-\sqrt{10a(6-a)+2(4+a)\sqrt{5a(4+a)}}}{4(a-1)}, \\
  c_4(a) &=& \frac{4+a+\sqrt{5a(4+a)}+\sqrt{10a(6-a)+2(4+a)\sqrt{5a(4+a)}}}{4(a-1)}.
\end{eqnarray*}

\begin{figure}[h!!]
    \centering
    \subfigure{
    \begin{tikzpicture}
    \begin{axis}[width=8cm,  axis equal image, scale only axis,  enlargelimits=false, axis on top]
      \addplot graphics[xmin=-55,xmax=85,ymin=-70,ymax=70] {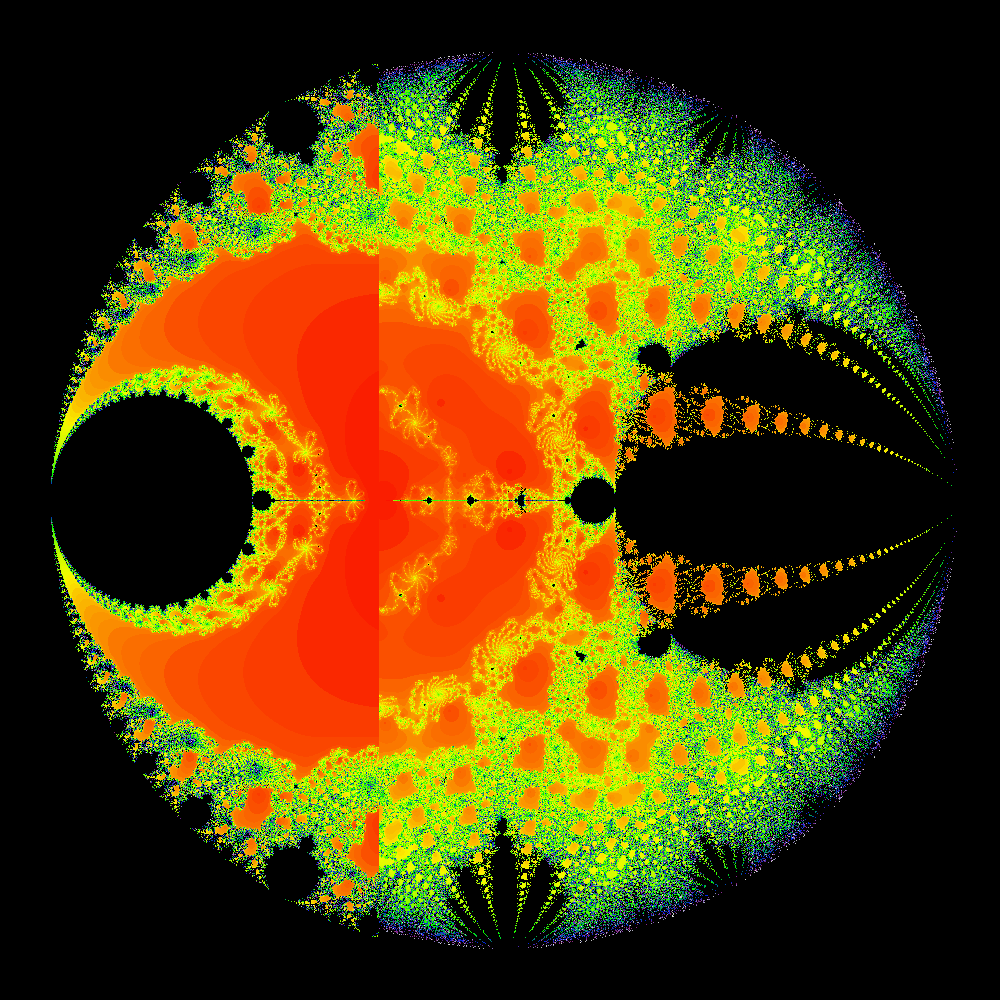};
    \end{axis}
  \end{tikzpicture}}
\subfigure{
	\begin{tikzpicture}
		\begin{axis}[width=8cm,  axis equal image,  scale only axis,  enlargelimits=false, axis on top]
			\addplot graphics[xmin=-55,xmax=85,ymin=-70,ymax=70] {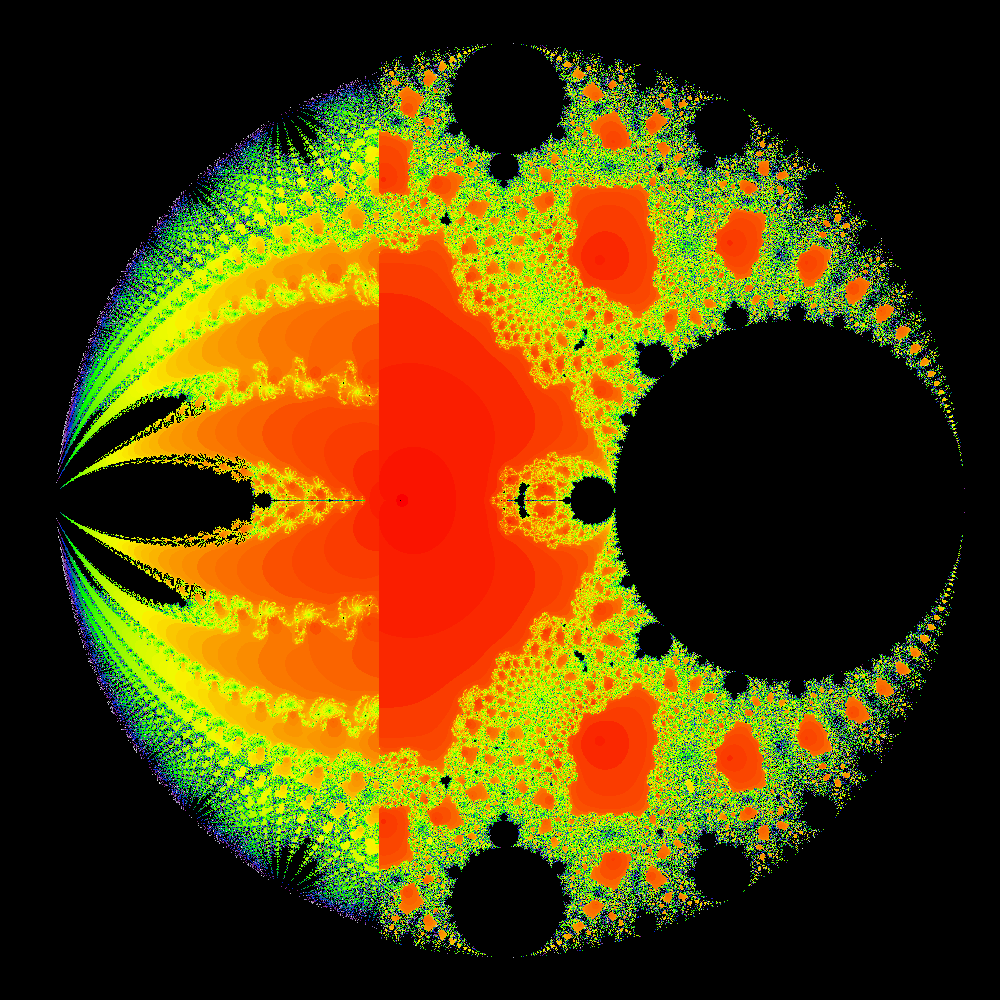};
		\end{axis}
\end{tikzpicture}}
\subfigure{
 \begin{tikzpicture}
	\begin{axis}[width=8cm,  axis equal image, scale only axis,  enlargelimits=false, axis on top]
		\addplot graphics[xmin=-55,xmax=85,ymin=-70,ymax=70] {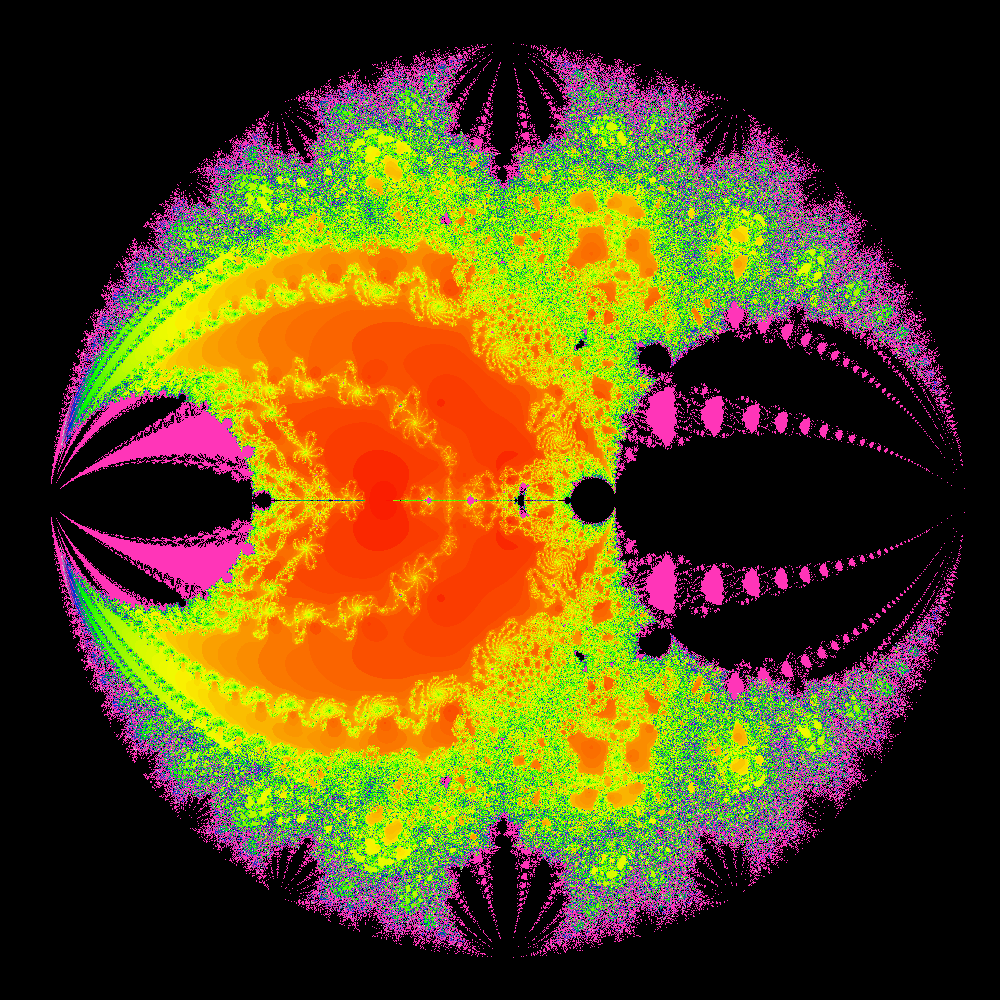};
	\end{axis}
\end{tikzpicture}
}
    \caption{\small{Upper figures show the parameter planes of the operator \eqref{op2} using the critical points $c_1(a)$ (left)  and $c_3(a)$ (right) separately. Lower figure shows the parameter plane obtained when using both critical points simultaneously.}}
    \label{c1c3OP2v2}
\end{figure}

As $c_2(a)=\frac{1}{c_1(a)}$ and $c_4(a)=\frac{1}{c_3(a)}$, it is enough to study the behaviour of $c_1$ and $c_3$. If we draw the parameter plane of each of them (see upper plots in Figure \ref{c1c3OP2v2}), we observe some inconsistencies   due to the indeterminacy generated by the square roots appearing in the critical points. 

These inconsistencies disappear when drawing a parameter plane that  takes into account both free critical points, as it can be observed in the lower plot in Figure \ref{c1c3OP2v2}. Recall that, when plotting the parameter plane using two critical points, black indicates that no critical orbit converges to the roots while pink indicates that only one critical orbit converges to the roots. If both critical orbits converge to the roots we use a scaling of colours depending on the slower time of convergence.

\begin{figure}[h!]
    \centering
    \subfigure[$a=-30$]{
    \begin{tikzpicture}
    \begin{axis}[width=6cm,  axis equal image, scale only axis,  enlargelimits=false, axis on top]
      \addplot graphics[xmin=-5,xmax=3,ymin=-4,ymax=4] {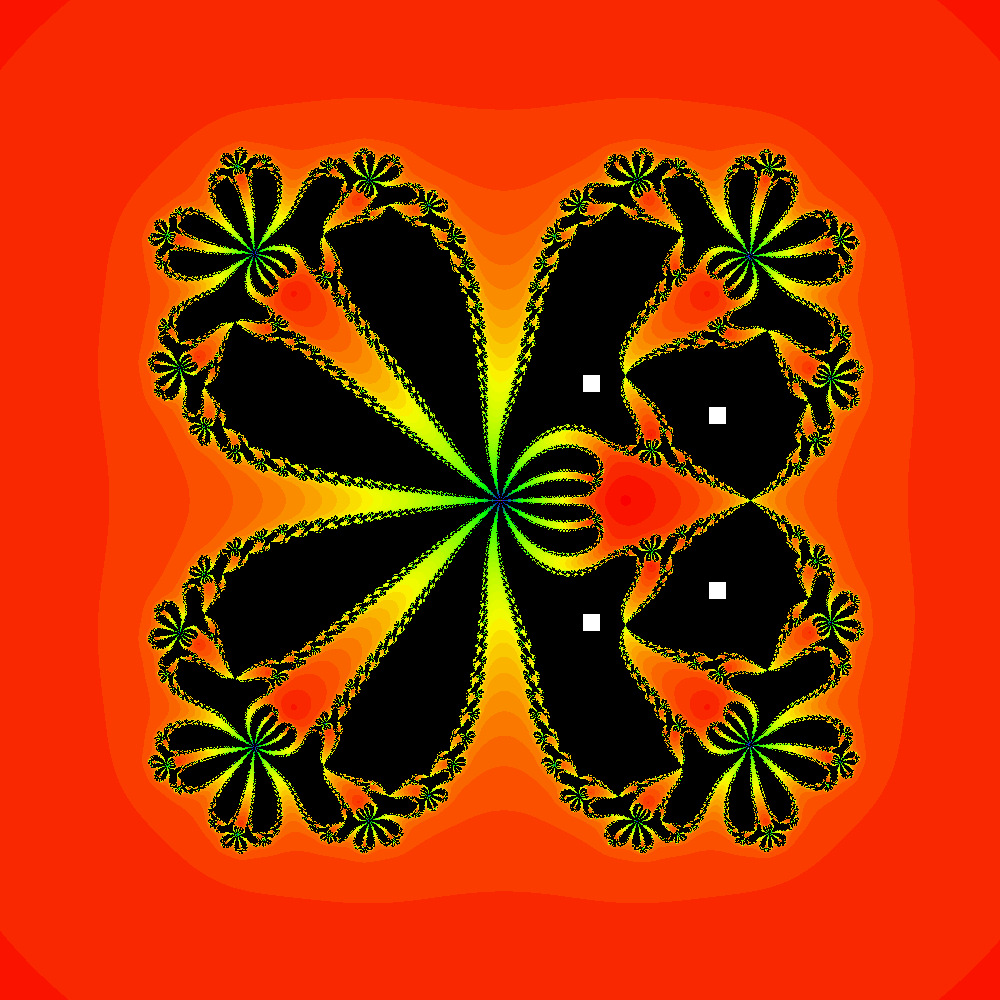};
    \end{axis}
  \end{tikzpicture}}
  \subfigure[$a=-30+8i$]{
	\begin{tikzpicture}
		\begin{axis}[width=6cm,  axis equal image, scale only axis,  enlargelimits=false, axis on top]
			\addplot graphics[xmin=-5,xmax=3,ymin=-4,ymax=4] {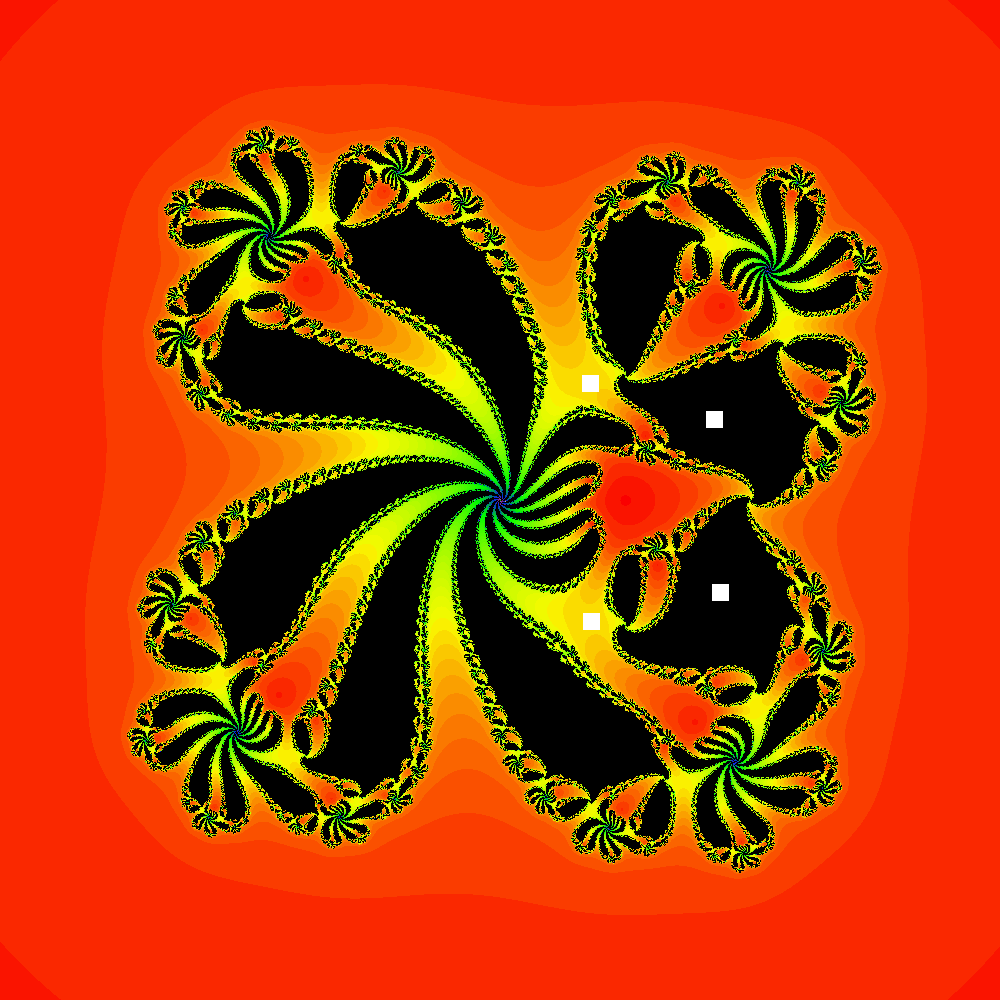};
		\end{axis}
\end{tikzpicture}}

   \subfigure[$a=-2$]{
    \begin{tikzpicture}
    \begin{axis}[width=6cm,  axis equal image, scale only axis,  enlargelimits=false, axis on top]
      \addplot graphics[xmin=-3,xmax=2,ymin=-2.5,ymax=2.5] {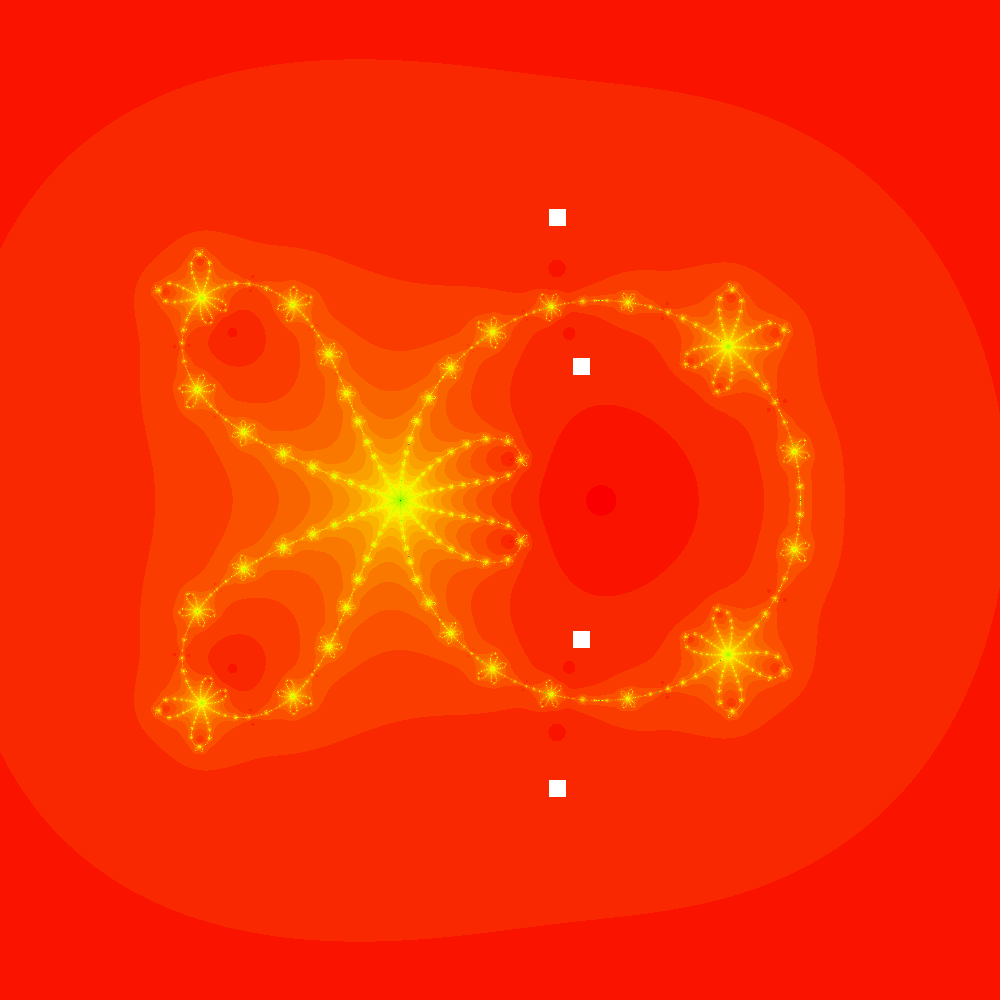};
    \end{axis}
  \end{tikzpicture}}
   \subfigure[$a=84$]{
    \begin{tikzpicture}
    \begin{axis}[width=6cm,  axis equal image, scale only axis,  enlargelimits=false, axis on top]
      \addplot graphics[xmin=-7,xmax=5,ymin=-6,ymax=6] {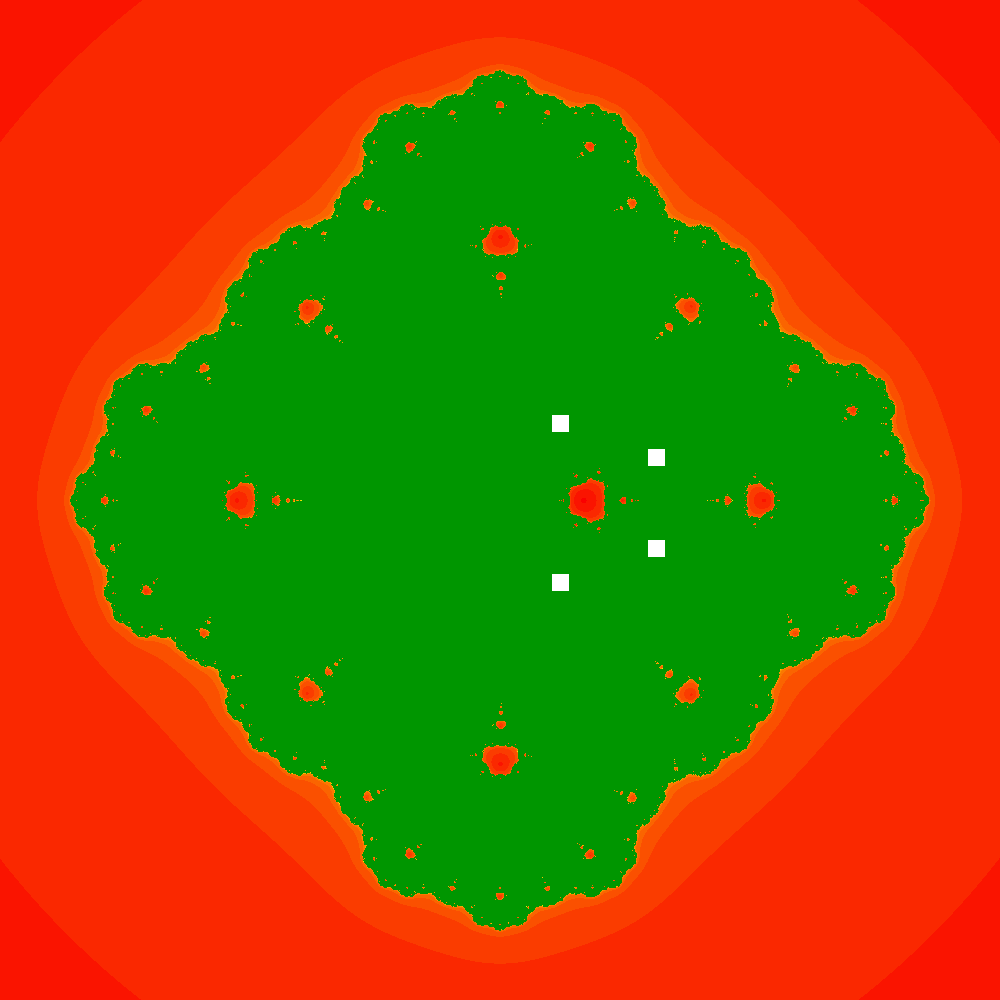};
    \end{axis}
  \end{tikzpicture}}
    \caption{\small{Dynamical planes of the operator \eqref{op2} for different values of the parameter $a$.}}
    \label{pdin_kim}
\end{figure}

In order to illustrate the different situations in the parameter plane, we finish this subsection by showing some dynamical planes of this operator (see Figure \ref{pdin_kim}). For these dynamical planes, we use the same scaling of colours used in the parameter planes to indicate convergence to the roots $z=0$ and $z=\infty$ (from red (fast convergence to the roots), to yellow, to pallid green, to blue and up to white (slow convergence), we use dark green if the point converges to $z=1$ (in case that $z=1$ is attracting), and we use black if the point does neither converge to the roots 0 and $\infty$ nor to the fixed point $z=1$. We also plot using white squares the location of the critical points. Notice that the family has four different free critical points (two modulo symmetry). Therefore, we will always have two critical points with the same dynamical behaviour.

 The parameters are chosen as follows. The value $a=-30$ belongs to a black parameter at the  bulb on the left. Since the parameter is black, no critical point converges to the roots. Indeed, the two critical points located more to the right  lie in the immediate basin of attraction of two different attracting fixed points (which are symmetric). The other two critical points belong the preimages of those basins of attraction. For $a=-30+8i$ the parameter is pink, so two critical points belong to the basins of attractions of the root. Indeed, this parameter belongs to the same bulb as $a=30$. In this case the two critical points to the left have moved out the basins of attraction of the fixed points, while the critical points to the right belong to immediate basins of attraction of the continuation of the attracting fixed points of $a=-30$.   
For $a=-2$ the parameter is red, so the four critical points are in the basin of attraction of the points $0$ and $\infty$, which correspond to the basins of attraction of the roots. Finally, the parameter $a=84$ lies on the unbounded black disk of parameters for which the point $z=1$ is attracting (and all free critical orbits converge to it).

\subsection{A multipoint variant of  Chebyshev's method}

In \cite{CCTV-varianteCheby}, the authors study the dynamics of a multipoint variant of Chebyshev's method. After applying it on quadratic polynomials they obtain the following operator:

\begin{equation}\label{op3}
    O_a(z)=\frac{z^3(2-4a+(5-8a+4a^2)z+(4-4a)z^2+z^3)}{(1+(4-4a)z+(5-8a+4a^2)z^2+(2-4a)z^3)}.
\end{equation}

The fixed points of  $O_a(z)$ are $z=0$, $z=\infty$ and the solutions of the equation:
\[
(z-1)(1 + (5 - 4a) z + 4(2 - 2 a + a^2)z^2 + (5 - 4 a)z^3 + z^4)=0.
\]
So, there are four strange fixed points in addition to $z=1$ that satisfy $z_1(a)=\frac{1}{z_2(a)}$ and $z_3(a)=\frac{1}{z_4(a)}$.

\begin{figure}[h]
\centering
\subfigure{
\includegraphics[scale=0.5]{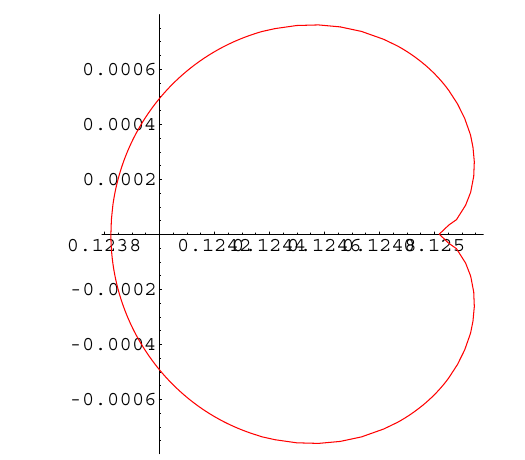}
}
\subfigure{
\includegraphics[scale=0.5]{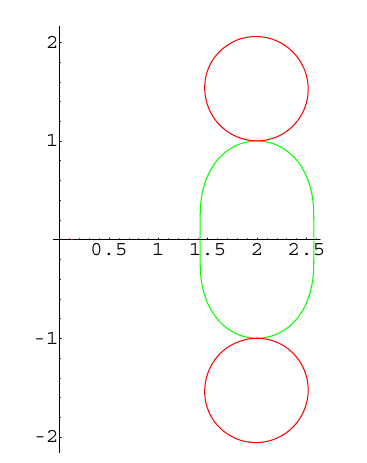}
}
\caption{Stability regions for $z=1$, $z_1$ and $z_2$ under the operator \eqref{op3}. }\label{est1OP3}
\end{figure}

The derivative of operator (\ref{op3}) is:
\begin{equation}\label{dop3}
    O^{\prime}_a(z)=\frac{-2z^2(1+(2-2a)z+z^2)P(z,a)}{(-1+(2a-1)z)^2(1+(3-2a)z+2z^2)^2},
\end{equation}
where
\[
P(z,a)=6a-3+(-12+22a-12a^2)z+(-18+32a-24a^2+8a^3)z^2+(-12+22a-12a^2)z^3+(6a-3)z^4.
\]
The fixed point $z=1$ is attractive inside the curve
\[
(55-24\alpha^3+3\alpha^4+22\beta^2+3\beta^4-8\alpha(13+3 \beta^2)+\alpha^2(74+6\beta^2)=0,
\]
being $a=\alpha+i\beta$. The point $z=-1$ is a preimage of $z=1$. We can also see in Figure \ref{est1OP3} the stability curves delimiting the regions where $z=1$ and the inverse pair $z_1$ and $z_2=1/z_1$ are attractive (coloured respectively in green and red). Notice that, by symmetry, the strange fixed points $z_1$ and $z_2$ are attractive for the same set of parameters.

From equation (\ref{dop3}), the points that satisfy  $1+(2-2a)z+z^2=0$ are preimages of $z=1$. So, the free critical points are the solutions  of $P(z,a)=0$. As in the previous subsection, the degree of the equation is reduced to the half by means of the change $z+\frac{1}{z}=x$. After solving the equation and undoing the change we have the four critical points:
\begin{eqnarray*}
  c_1(a) &=& \frac{6-11a+6a^2-a\sqrt{1+36a-12a^2}-\sqrt{2a \left( 6+25a-48a^2+12a^3-(6-11a+6a^2)\sqrt{1+36a-12a^2} \right) }}{6(-1+2a)}, \\
  c_2(a) &=&  \frac{6-11a+6a^2-a\sqrt{1+36a-12a^2}+\sqrt{2a  \left(6+25a-48a^2+12a^3-(6-11a+6a^2)\sqrt{1+36a-12a^2} \right) }}{6(-1+2a)}, \\
  c_3(a) &=&  \frac{6-11a+6a^2+a\sqrt{1+36a-12a^2}-\sqrt{2a  \left(6+25a-48a^2+12a^3+(6-11a+6a^2)\sqrt{1+36a-12a^2} \right) }}{6(-1+2a)}, \\
  c_4(a) &=&  \frac{6-11a+6a^2+a\sqrt{1+36a-12a^2}+\sqrt{2a  \left(6+25a-48a^2+12a^3+(6-11a+6a^2)\sqrt{1+36a-12a^2} \right) }}{6(-1+2a)}.
\end{eqnarray*}

\begin{figure}[h!!]
    \centering
    \subfigure{
    \begin{tikzpicture}
    \begin{axis}[width=9cm,  axis equal image, scale only axis,  enlargelimits=false, axis on top]
      \addplot graphics[xmin=-1.5,xmax=3.9,ymin=-4.5,ymax=4.5] {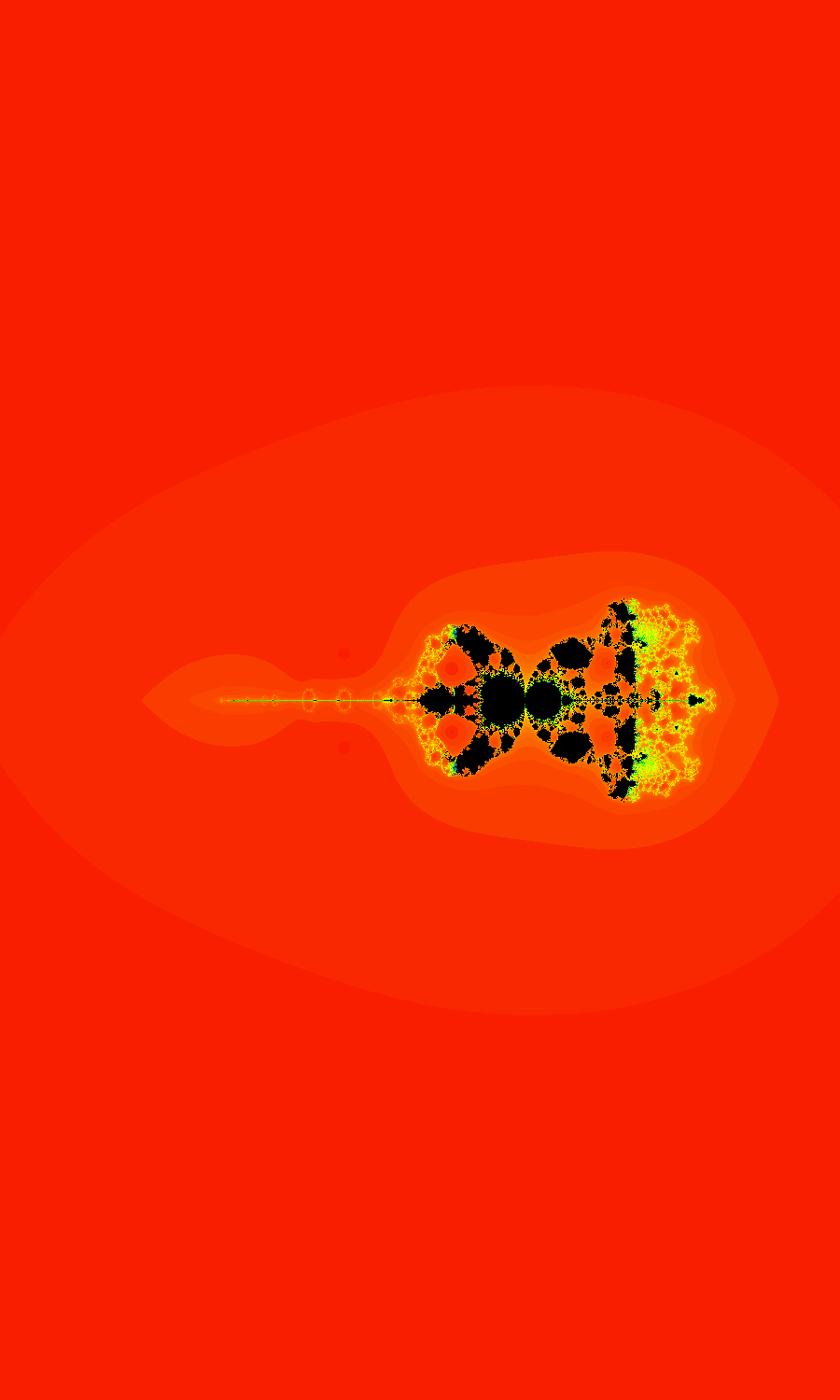};
    \end{axis}
  \end{tikzpicture}}
\subfigure{
	\begin{tikzpicture}
		\begin{axis}[width=9cm,  axis equal image, scale only axis,  enlargelimits=false, axis on top]
			\addplot graphics[xmin=-1.5,xmax=3.9,ymin=-4.5,ymax=4.5] {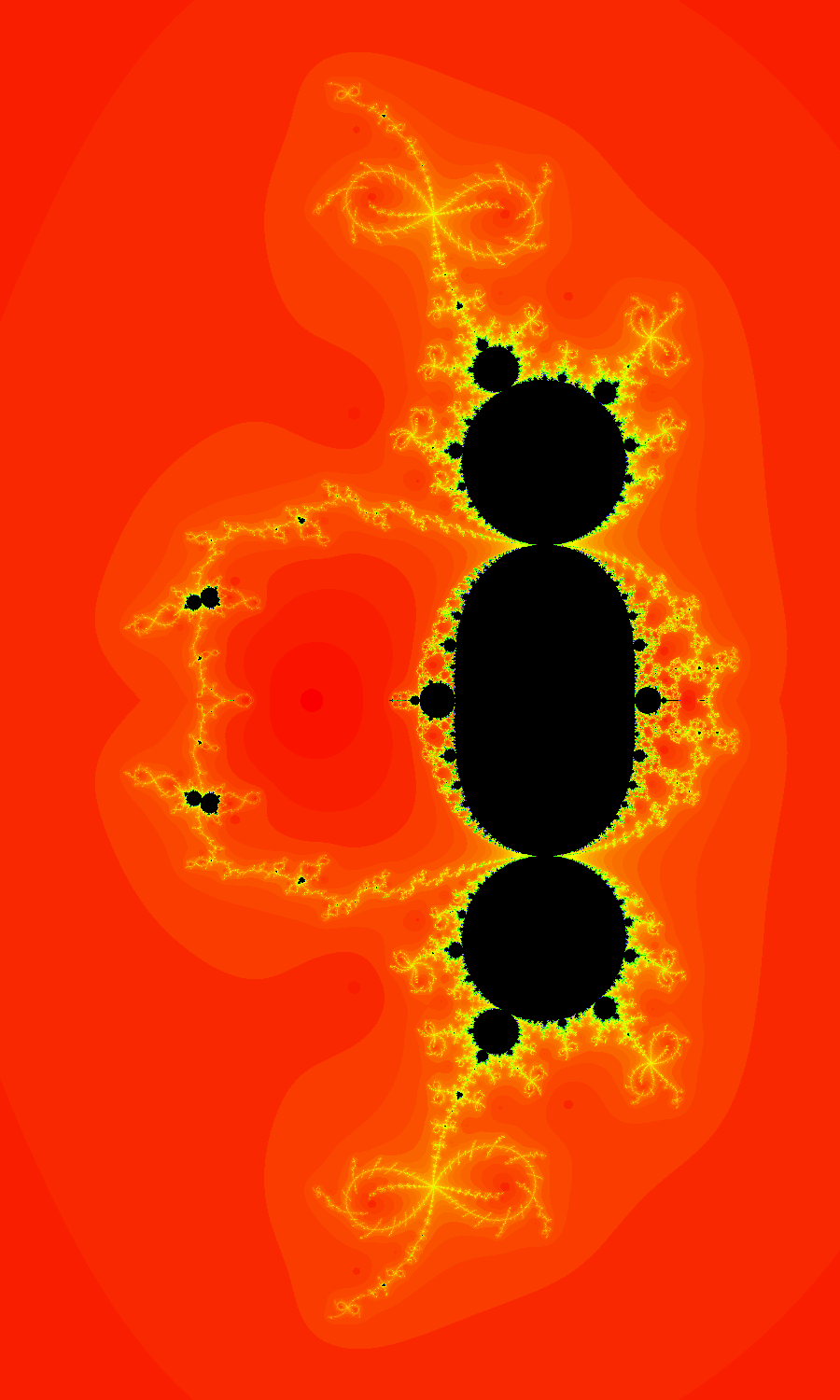};
		\end{axis}
\end{tikzpicture}}
\subfigure{
	\begin{tikzpicture}
		\begin{axis}[width=9cm,  axis equal image, scale only axis,  enlargelimits=false, axis on top]
			\addplot graphics[xmin=-1.5,xmax=3.9,ymin=-4.5,ymax=4.5] {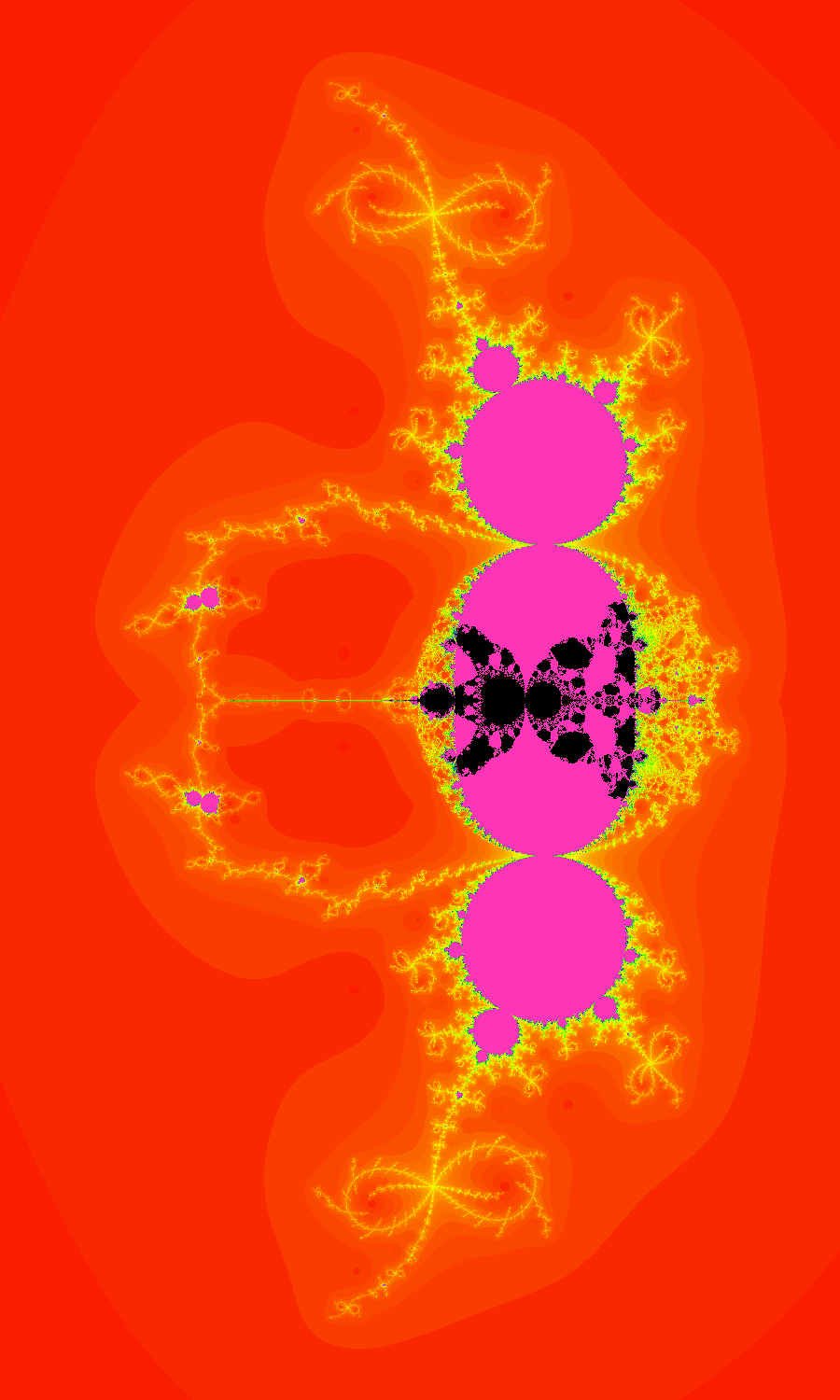};
		\end{axis}
\end{tikzpicture}}
    \caption{\small{Parameter planes of the operator \eqref{op3} obtained using the critical points $c_1(a)$ (left), the critical point $c_3(a)$ (centre), and both critical points simultaneously (right).}}
    \label{fig:param3}
\end{figure}

As $c_2(a)=\frac{1}{c_1(a)}$ and $c_4(a)=\frac{1}{c_3(a)}$,  it is enough to study the behaviour of $c_1$ and $c_3$. In Figure \ref{fig:param3} left and centre we show the parameter planes obtained using $c_1$ and $c_3$, respectively. In this case we do not observe any incongruence coming from changes of determination of the critical points. However, by obtaining a parameter plane using both critical points (see Figure \ref{fig:param3} right), we do obtain a better understanding of the bifurcation set.

Recall that in Figure~\ref{fig:param3} right, colour black indicates that no free critical orbit converges to the roots while pink indicates that only one free critical orbit converges to the root. Therefore, if the parameter is black there may be up to four basins of attraction other than the roots (two up to symmetry). It is then relevant to know if the critical points $c_1$ and $c_3$ converge to two different  attracting cycles (not related by symmetry), since those parameters would be particularly inconvenient.

\begin{figure}[h!!]
	\centering
	\subfigure{
		\begin{tikzpicture}
			\begin{axis}[width=9cm, axis equal image,  scale only axis,  enlargelimits=false, axis on top]
				\addplot graphics[xmin=-1.5,xmax=3.9,ymin=-4.5,ymax=4.5] {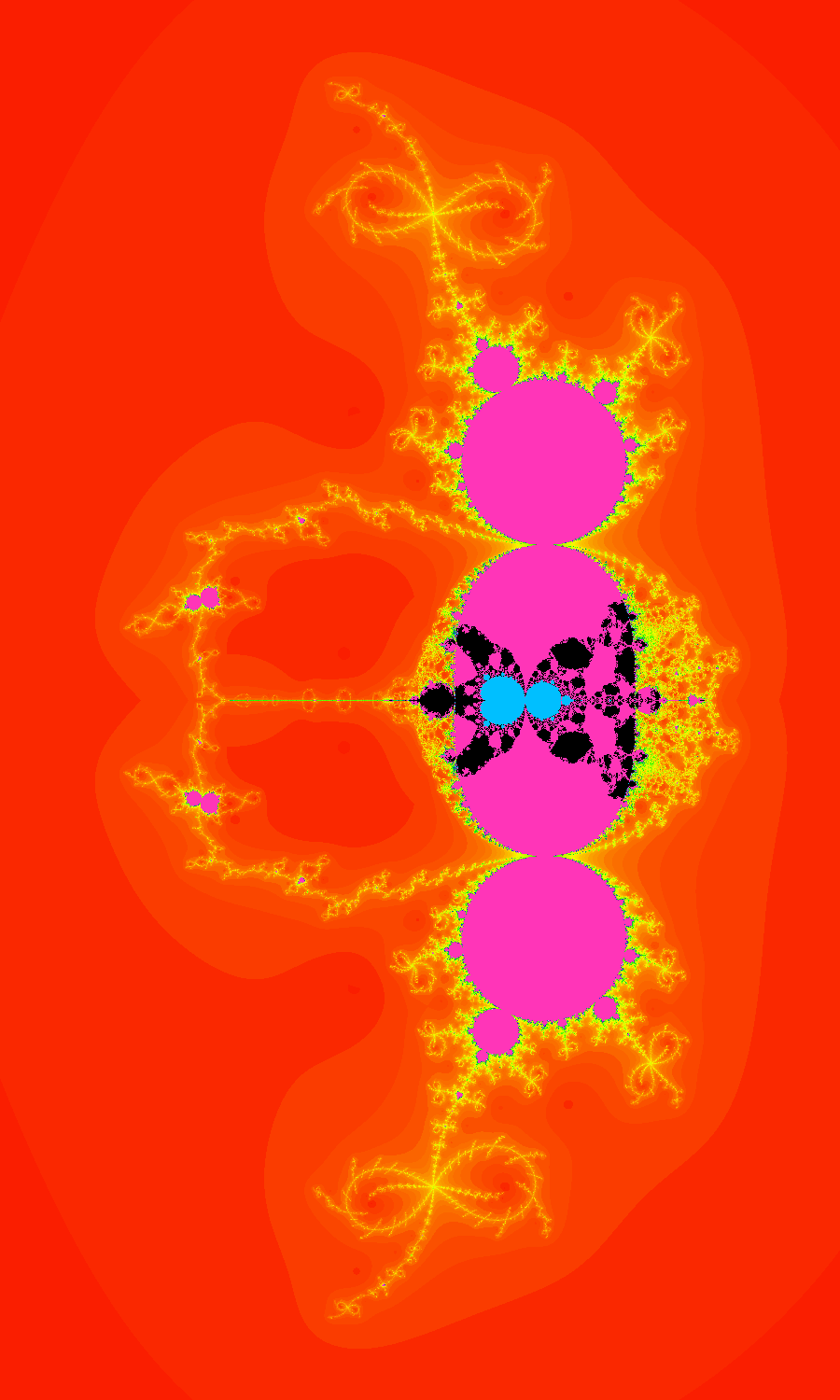};
			\end{axis}
	\end{tikzpicture}}
	\subfigure{
		\begin{tikzpicture}
			\begin{axis}[width=9cm,  axis equal image, scale only axis,  enlargelimits=false, axis on top]
				\addplot graphics[xmin=1.25,xmax=2.75,ymin=-0.75,ymax=0.75] {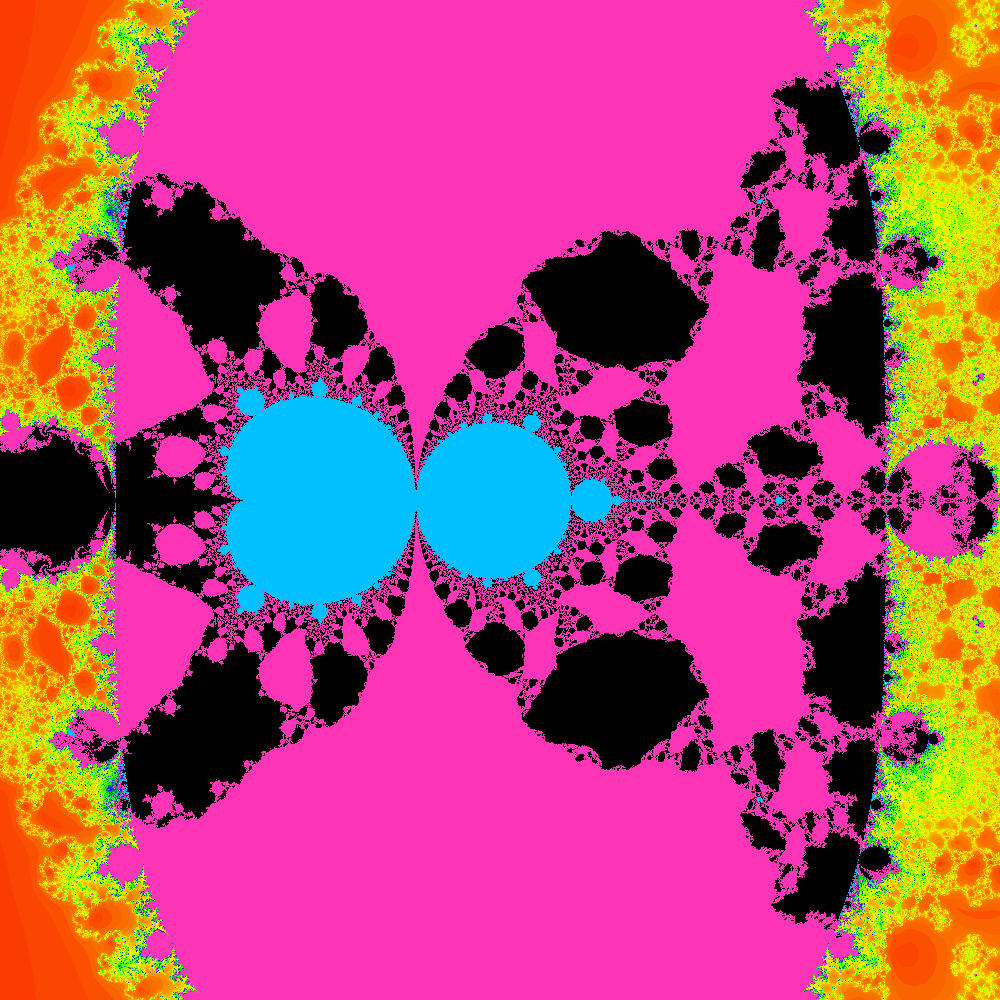};
			\end{axis}
	\end{tikzpicture}}
	\caption{\small{Parameter plane of the operator \eqref{op3} obtained using the critical points $c_1$ and $c_3$. We use colour pink if only one free critical orbit converges to the roots, blue if the orbits od $c_1$ and $c_3$ converge to two different non-symmetric attracting cycles, and black if both free critical orbits converge to symmetric attracting cycles. }}
	\label{fig:param3-2col}
\end{figure}

The program we use can be modified to detect if $c_1$ and $c_3$ converge to different cycles. If neither $c_1$ nor $c_3$ converge to any of the roots after the allowed maximal number of iterates, we check if they converge to the same cycle  (or a symmetric one). To do so we first iterate $c_1$ up to $10000$ times to ``refine'' the set at which it converges. Afterwards we store in two different vectors the 100 first iterates of $c_1$ and $c_2=1/c_1$. Next, we iterate $c_3$ up to  10000 times, obtaining a point $z_3$ of its orbit. Finally, we verify if $z_3$ coincides with any of the 100 iterates   that we stored of the orbits of $c_1$ and $c_2$. Notice that using this procedure we can only detect that they converge to the same cycle if its period is smaller than 100. The result of this procedure can be observed in Figure~\ref{fig:param3-2col}. As before, pink indicates that only one critical value converges to the roots. Black indicates that both $c_1$ and $c_3$ converge to the same cycle (up to symmetry). Blue indicates that $c_1$ and $c_3$ converge to different cycles. Parameters corresponding to black points, for which there is more than one free critical orbit within the same basin of attraction, are usually called captures parameters. On the other hand, if $c_1$ and $c_3$ converge to different cycles we call the parameter disjoint. We would like to point out that disjoint parameters can usually be recognized without having to use this colouring. Indeed, disjoint parameters usually lead to Mandelbrot-like structures in the parameter plane which are easily recognizable (see Figure~\ref{fig:param3-2col}). We would like to point out that the black components in Figure~\ref{c1c3OP2v2} are mostly capture components (disjoint hyperbolic components of Operator \ref{op2} are very small).

\begin{figure}[h!]
	\centering
	\subfigure[$a=1.45$]{
		\begin{tikzpicture}
			\begin{axis}[width=6cm,  axis equal image, scale only axis,  enlargelimits=false, axis on top]
				\addplot graphics[xmin=-2,xmax=3,ymin=-2.5,ymax=2.5] {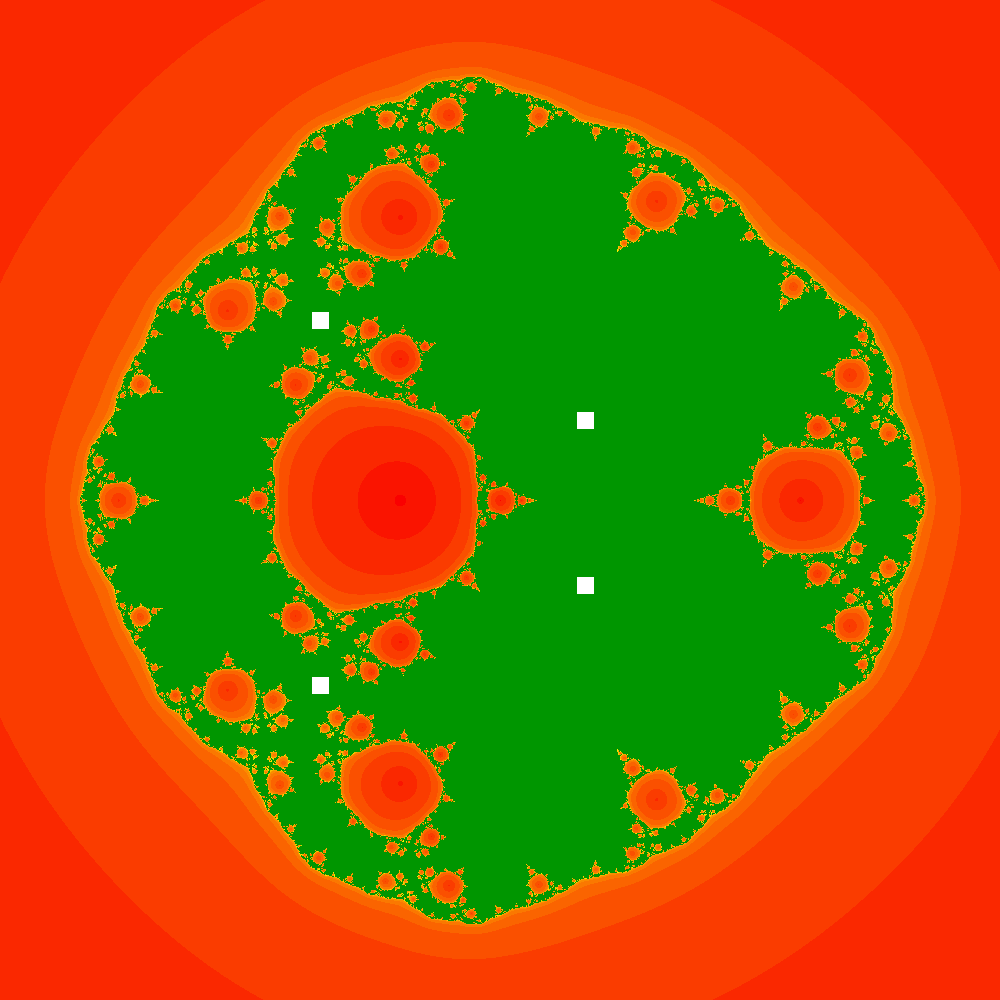};
			\end{axis}
	\end{tikzpicture}}
	\subfigure[$a=1.7$]{
		\begin{tikzpicture}
			\begin{axis}[width=6cm,  axis equal image, scale only axis,  enlargelimits=false, axis on top]
				\addplot graphics[xmin=-2,xmax=3.5,ymin=-2.75,ymax=2.75] {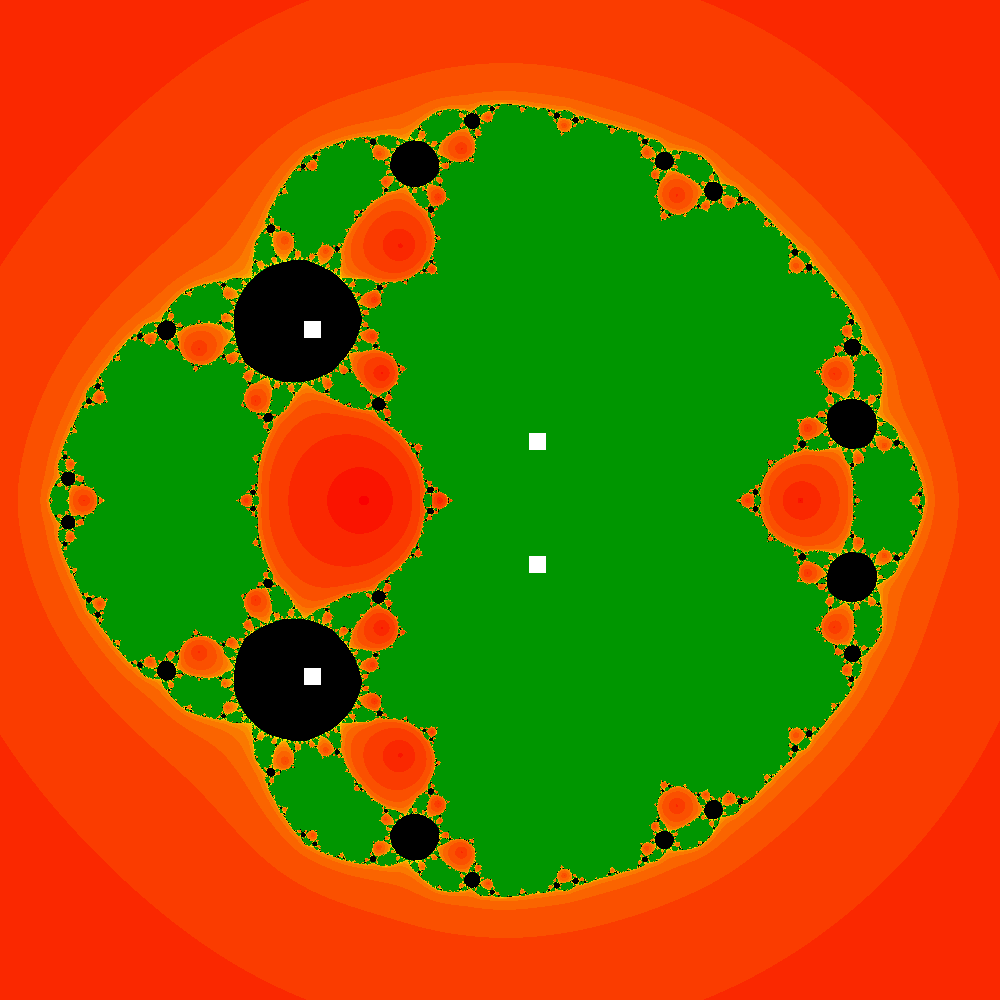};
			\end{axis}
	\end{tikzpicture}}
	\subfigure[$a=2+1.5i$]{
		\begin{tikzpicture}
			\begin{axis}[width=6cm, axis equal image, scale only axis,  enlargelimits=false, axis on top]
				\addplot graphics[xmin=-3,xmax=5,ymin=-2,ymax=6] {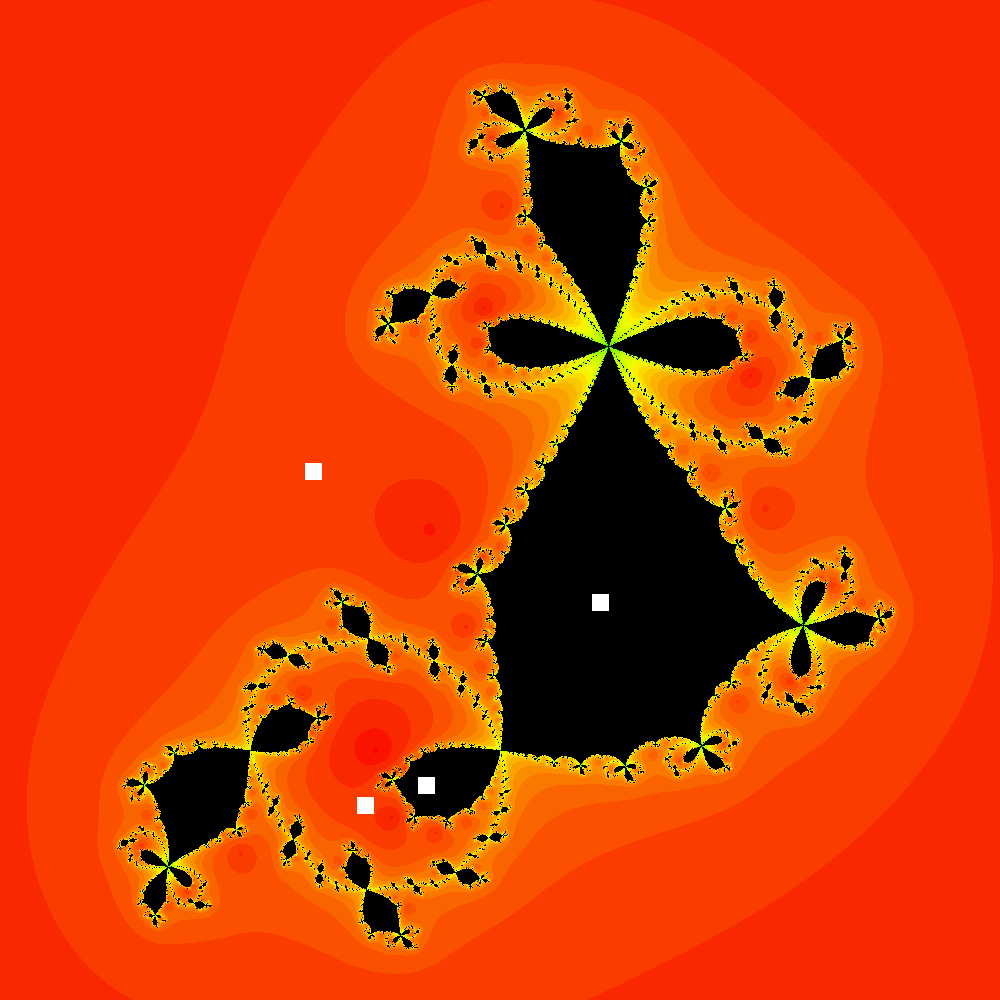};
			\end{axis}
	\end{tikzpicture}}
	\subfigure[$a=1.9+0.5i$]{
		\begin{tikzpicture}
			\begin{axis}[width=6cm,  axis equal image, scale only axis,  enlargelimits=false, axis on top]
				\addplot graphics[xmin=-5.5,xmax=3.5,ymin=-3,ymax=6] {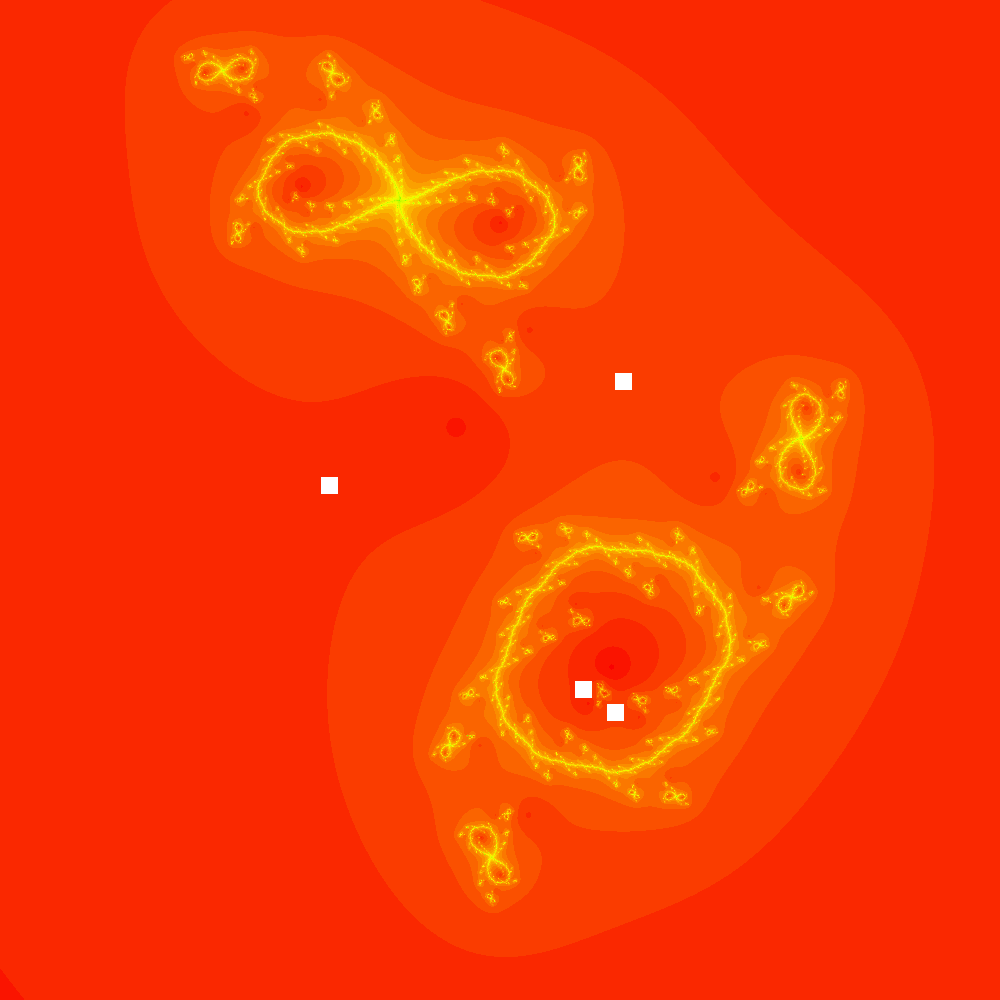};
			\end{axis}
	\end{tikzpicture}}
	\caption{\small{Dynamical planes of the operator \eqref{op3} for different values of the parameter $a$.} }
	\label{pdin_cheby}
\end{figure}

We end this subsection with some dynamical planes that illustrate the different colours in the parameter plane of Operator \ref{op3}  (see Figure \ref{pdin_cheby}). The colours are as in the dynamical planes in \S~\ref{subsec:op2}.  The value $a=1.5$ is located in a black region inside the main bulb of Figure~\ref{fig:param3-2col}. All four critical points (two up to symmetry) belong to the basin of attraction of $z=1$ (in green). The parameter $a=1.7$ lies in the biggest blue cardioid of Figure~\ref{fig:param3-2col}. The basin of attraction of $z=1$ (in green) contains two critical points (one up to symmetry). The other two free critical points lie in the basin of attraction of an attracting cycle of period two (in black). The parameter  $a=2+1.5i$ belongs to a pink bulb where the fixed points $z_1$ and $z_2$ are attracting  (compare Figure \ref{est1OP3} (right)). In black we observe the basins of attraction of these two points, each of them containing a free critical point. The other two free critical points belong to the basins of attraction of the roots $z=0$ and $z=\infty$. Finally, the parameter $a=2i$ (which appears as a red parameter) is chosen so that all critical points belong to the basins of attraction of the roots.

 \subsection{The Ermakov-Kalitkin family} \label{seccion_Ermakov-Kalitkin family}

The next example of application of the program corresponds to a family with two free critical orbits (up to symmetry) from which only one orbit is actually free. In \cite{CTV-ermakov} the authors study the dynamical behaviour of the   family of Ermakov-Kalitkin type methods applied to quadratic polynomials. In this case, the strange fixed point $z=-1$ is a parabolic point of multiplicity three and it is located on the boundary of two parabolic basins. It follows that each of these parabolic basins must contain a critical point.

After applying the method on quadratic polynomials the following operator is obtained:

\begin{equation}\label{eq:op4}
    O_a(z)=\frac{z^3(2(a-1)+a^2+4(a-1)z+2(a-1)z^2)}{2(a-1)+4(a-1)z+(a^2+2a-2)z^2},
\end{equation}
whose derivative is:
\begin{equation}\label{dop4}
    O^{\prime}_a(z)=\frac{z^2 P(a,z)} {(2(a-1)+4(a-1)z+(a^2+2a-2)z^2)^2)},
\end{equation}
being
\begin{eqnarray*}\label{dopp4}
   P(a,z)&=& 6(a-1)(-2+2a+a^2)+8(a-1)(-6+6a+a^2)z+(72-144a+68a^2+4a^3+a^4)z^2 +   \\
  &+ &  8(a-1)(-6+6a+a^2)z^3+6(a-1)(-2+2a+a^2)z^4.
\end{eqnarray*}

In this case, the strange fixed points are only $z=1$ and $z=-1$, since $O_a(z)=z$ implies:
\[
2(-1+ a)(-1 + z)(1 + z)^3=0.
\]

\begin{figure}[h]
\centering
\subfigure{
\includegraphics[scale=0.5]{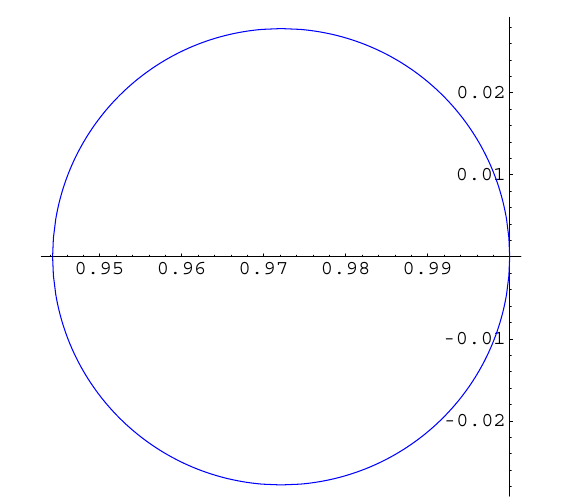}
}
\subfigure{
\includegraphics[scale=0.8]{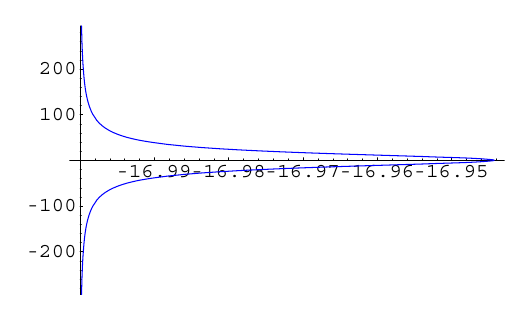}
}
\caption{Stability regions for $z=1$.}\label{estErmakov}
\end{figure}

The point $z=1$ is attractive inside the curve defined by
\[
16-32\alpha+15\alpha^2+\alpha^3+17\beta^2+\alpha \beta^2=0,
\]
being $a=\alpha+i\beta$ (see Figure \ref{estErmakov}).

The fixed point $z=-1$ is a parabolic point with multiplicity 3, as it is a triple solution of $O_a(z)=z$. Therefore, $z=-1$ lies on the boundary of two attractive parabolic basins, each of them containing a critical point (see \cite[\S 10]{Milnor}, for instance).

\begin{figure}[hbt!]
	\centering
	\subfigure{
		\begin{tikzpicture}
			\begin{axis}[width=8cm,  axis equal image, scale only axis,  enlargelimits=false, axis on top]
				\addplot graphics[xmin=-27,xmax=13,ymin=-20,ymax=20] {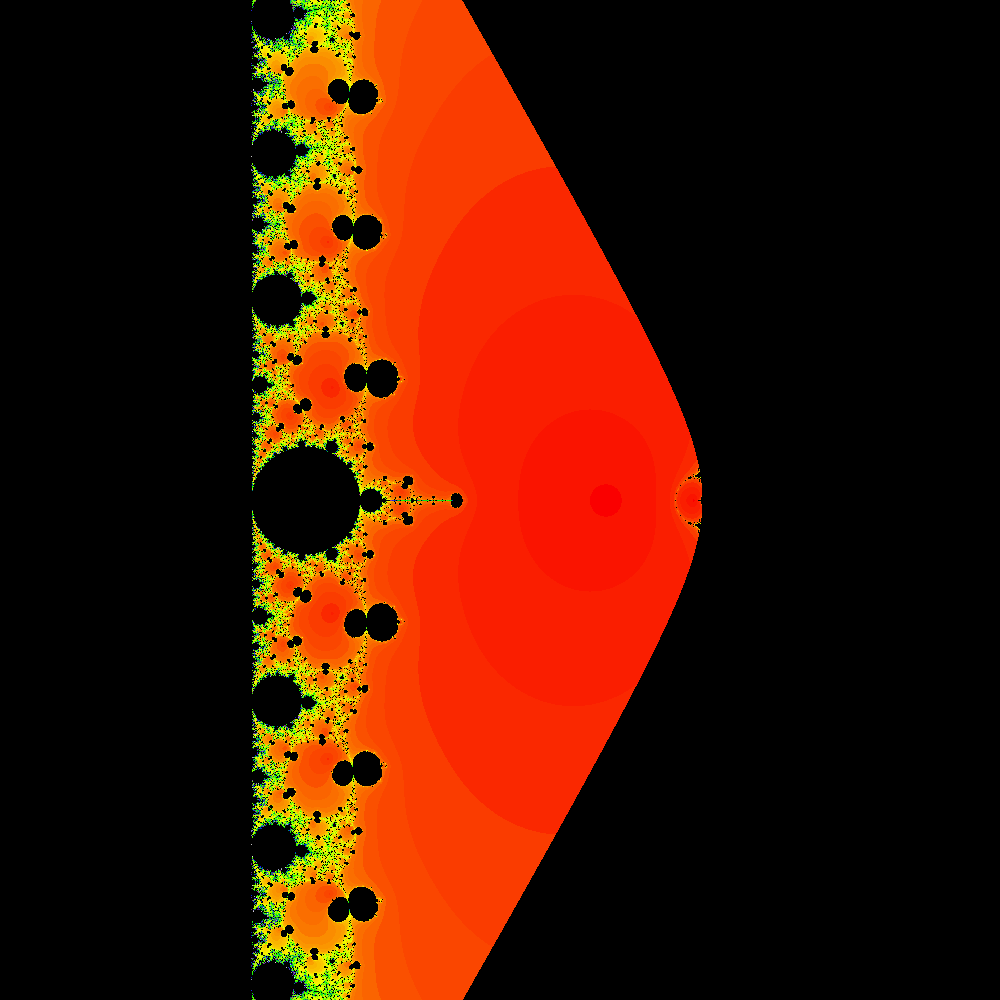};
			\end{axis}
	\end{tikzpicture}}
	\subfigure{
		\begin{tikzpicture}
			\begin{axis}[width=8cm, axis equal image, scale only axis,  enlargelimits=false, axis on top]
				\addplot graphics[xmin=-27,xmax=13,ymin=-20,ymax=20] {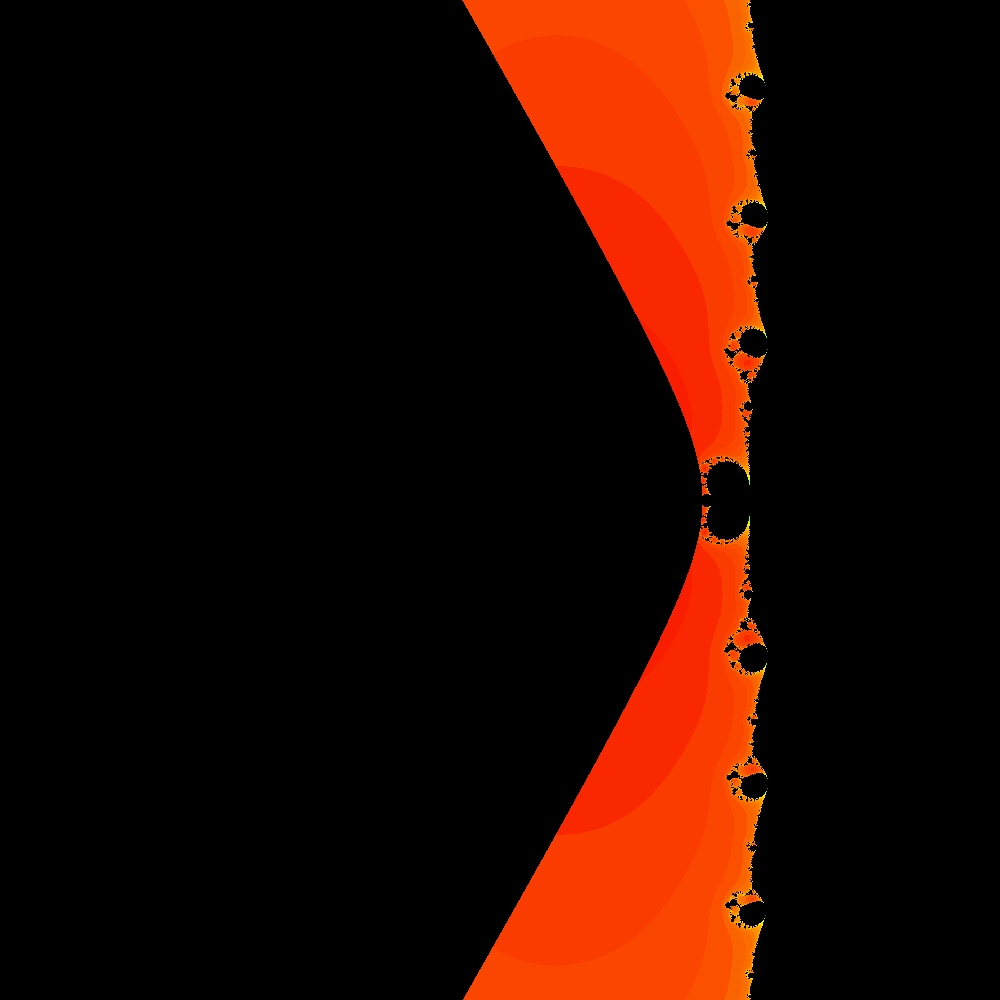};
			\end{axis}
	\end{tikzpicture}}
	
	\subfigure{
		\begin{tikzpicture}
			\begin{axis}[width=8cm, axis equal image, scale only axis,  enlargelimits=false, axis on top]
				\addplot graphics[xmin=-27,xmax=13,ymin=-20,ymax=20] {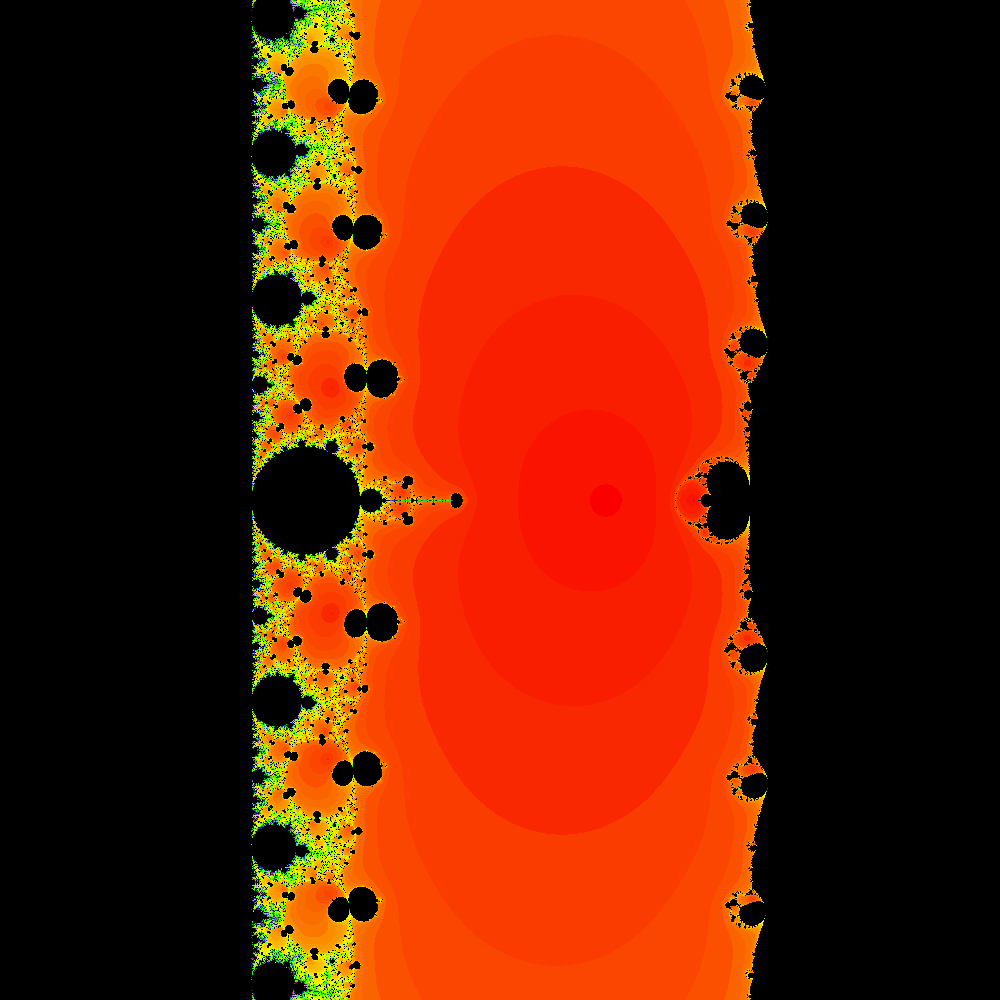};
			\end{axis}
	\end{tikzpicture}}
	\caption{\small{Parameter planes of the operator \eqref{eq:op4} obtained using the critical points $c_1(a)$ (upper left) and $c_3(a)$ (upper right) separately and using both critical points simultaneously (lower).}}
	\label{fig:paramO4}
\end{figure}

As in the previous examples, we can obtain the four roots of the polynomial $P(a,z)$ in (\ref{dopp4}). They correspond to the four free critical points:
\begin{eqnarray*}
  c_1(a) &=& \frac{4(1-a)(-6+6 a+a^2)-a^2\sqrt{2(1-a)(26-26 a+3 a^2)}}{12(a-1)(-2+2 a+a^2)}- \\
             &-&\frac{a \sqrt{2(1-a)\left( 192-384 a+154 a^2+38 a^3+3 a^4-4(-6+6 a+a^2)\sqrt{2(1-a)(26-26 a+3 a^2)}\right) }}{12(a-1)(-2+2 a+a^2)}, 
 \end{eqnarray*}
\begin{eqnarray*}
  c_2(a) &=& \frac{4(1-a)(-6+6a+a^2)-a^2\sqrt{2(1-a)(26-26 a+3 a^2)}}{12(a-1)(-2+2  a+a^2)}+ \\
             &+&\frac{a \sqrt{2(1-a)\left(192-384 a+154 a^2+38 a^3+3 a^4-4(-6+6 a+a^2)\sqrt{2(1-a)(26-26 a+3 a^2)}\right) }}{12(a-1)(-2+2 a+a^2)}, 
   \end{eqnarray*}
  \begin{eqnarray*}
  c_3(a) &=& \frac{4(1-a)(-6+6 a+a^2)+a^2\sqrt{2(1-a)(26-26 a+3 a^2)}}{12(a-1)(-2+ a+a^2)}- \\
             &-&\frac{a \sqrt{2(1-a)\left(192-384 a+154 a^2+38 a^3+3 a^4+4 (-6+6 a+a^2) \sqrt{2(1-a)(26-26 a+3 a^2)}\right) }}{12(a-1)(-2+2 a+a^2)}, 
   \end{eqnarray*}
\begin{eqnarray*}
  c_4(a) &=& \frac{4(1-a)(-6+6 a+a^2)+a^2 \sqrt{2(1-a)(26-26 a+3 a^2)}}{12(a-1)(-2+2 a+a^2)}+ \\
             &+&\frac{a \sqrt{2(1-a)\left(192-384 a+154 a^2+38 a^3+3 a^4+4(-6+6 a+a^2) \sqrt{2(1-a)(26-26 a+3 a^2)}\right) }}{12(a-1)(-2+2 a+a^2)}.
\end{eqnarray*}

As $c_2(a)=\frac{1}{c_1(a)}$ and $c_4(a)=\frac{1}{c_3(a)}$, it is enough to study the behaviour of $c_1$ and $c_3$. If we plot the parameter plane of each of them (see Figure \ref{fig:paramO4} upper) we observe the inconsistencies  due to the interaction between the two critical points.

Since $z=-1$  is a parabolic point of multiplicity three and therefore, it is located at the boundary two parabolic basins. One critical point must lie in each of these parabolic basins. However this may take only one free critical orbit up to symmetry, since it is enough that one of the critical points and its inverse lie in each of these basins. It follows that at most one of the critical points (up to symmetry) may lead to new stable dynamics (other than the basins of attraction of the roots and $z=-1$). 

The parameter plane obtained using the two free critical points $c_1$ and $c_3$ can be observed in Figure~\ref{fig:paramO4} (lower). However, since only one of the orbits may actually be free (one of the critical orbits has to converge to the parabolic point $z=-1$), it is convenient to do a small modification to the program. In the modified version if no critical point converges to the roots we use black and we use the scaling of colours that we usually use when all critical orbits converge to the roots when one critical orbit converges to the roots (that is the maximum amount critical orbits that can escape to the roots). If we did not do this modification then regions with scaling of red would appear as pink (which would still be correct but would make more difficult to recognise the thinner regions of bifurcation parameters).

\begin{figure}[h!]
	\centering
	\subfigure[$a=-20$]{
		\begin{tikzpicture}
			\begin{axis}[width=6cm,  axis equal image, scale only axis,  enlargelimits=false, axis on top]
				\addplot graphics[xmin=-6,xmax=4,ymin=-5,ymax=5] {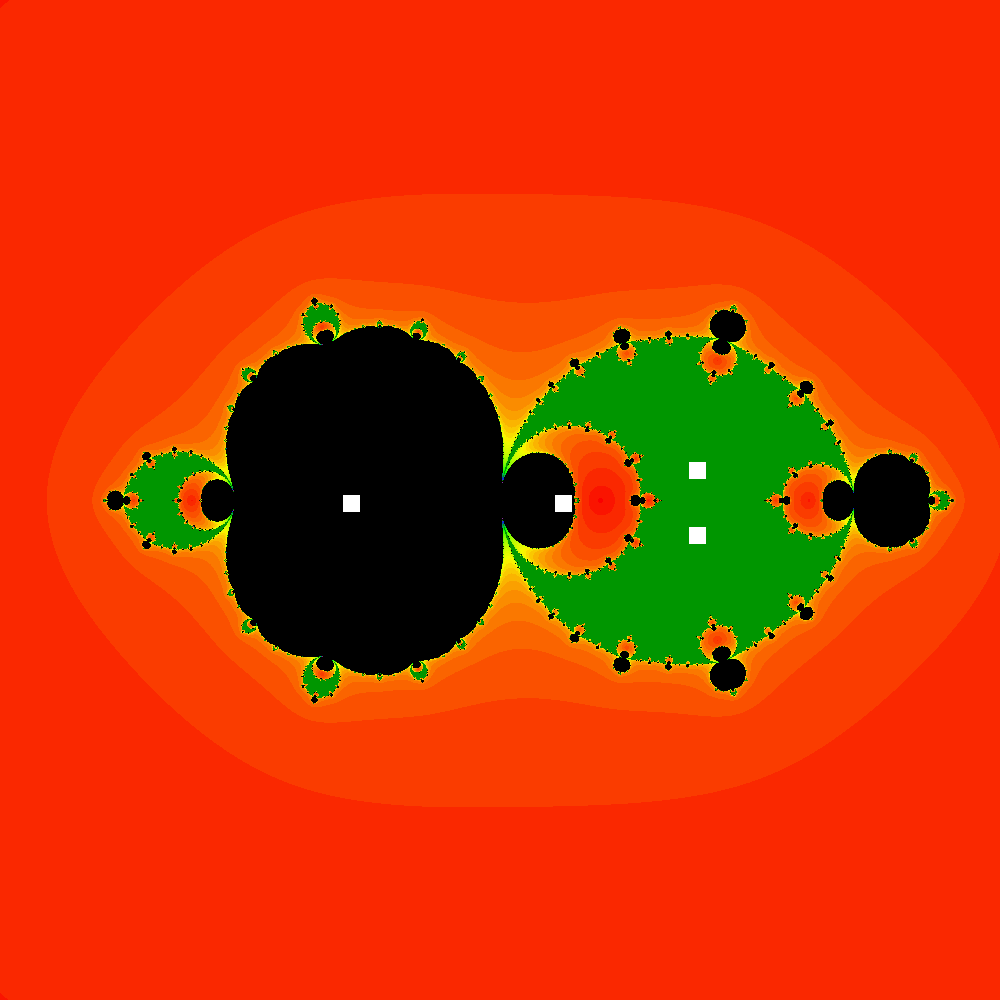};
			\end{axis}
	\end{tikzpicture}}
	\subfigure[$a=0.97$]{
		\begin{tikzpicture}
			\begin{axis}[width=6cm,  axis equal image, scale only axis,  enlargelimits=false, axis on top]
				\addplot graphics[xmin=-6.5,xmax=5.5,ymin=-6,ymax=6] {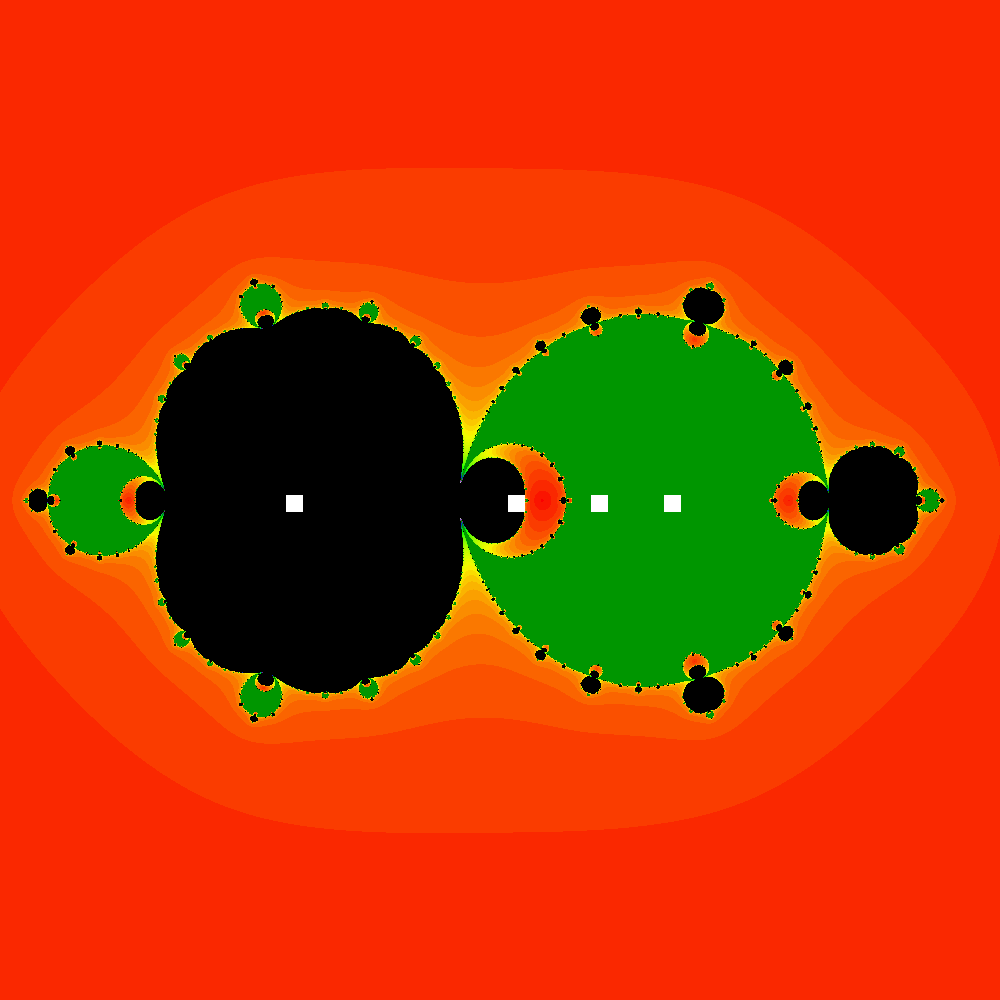};
			\end{axis}
	\end{tikzpicture}}
	
	\subfigure[$a=-7$]{
		\begin{tikzpicture}
			\begin{axis}[width=6cm,  axis equal image, scale only axis,  enlargelimits=false, axis on top]
				\addplot graphics[xmin=-4,xmax=2,ymin=-3,ymax=3] {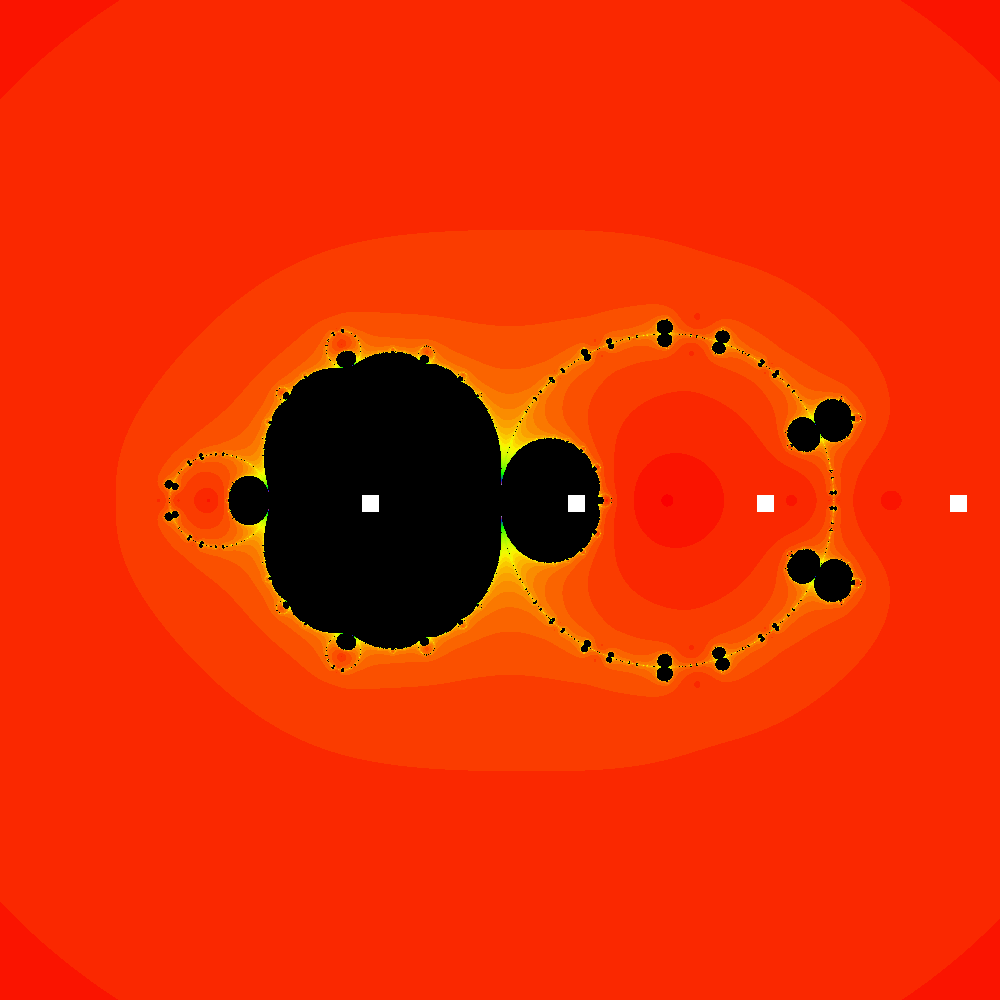};
			\end{axis}
	\end{tikzpicture}}
	\subfigure[$a=9$]{
		\begin{tikzpicture}
			\begin{axis}[width=6cm, axis equal image, scale only axis,  enlargelimits=false, axis on top]
				\addplot graphics[xmin=-4,xmax=3,ymin=-3.5,ymax=3.5] {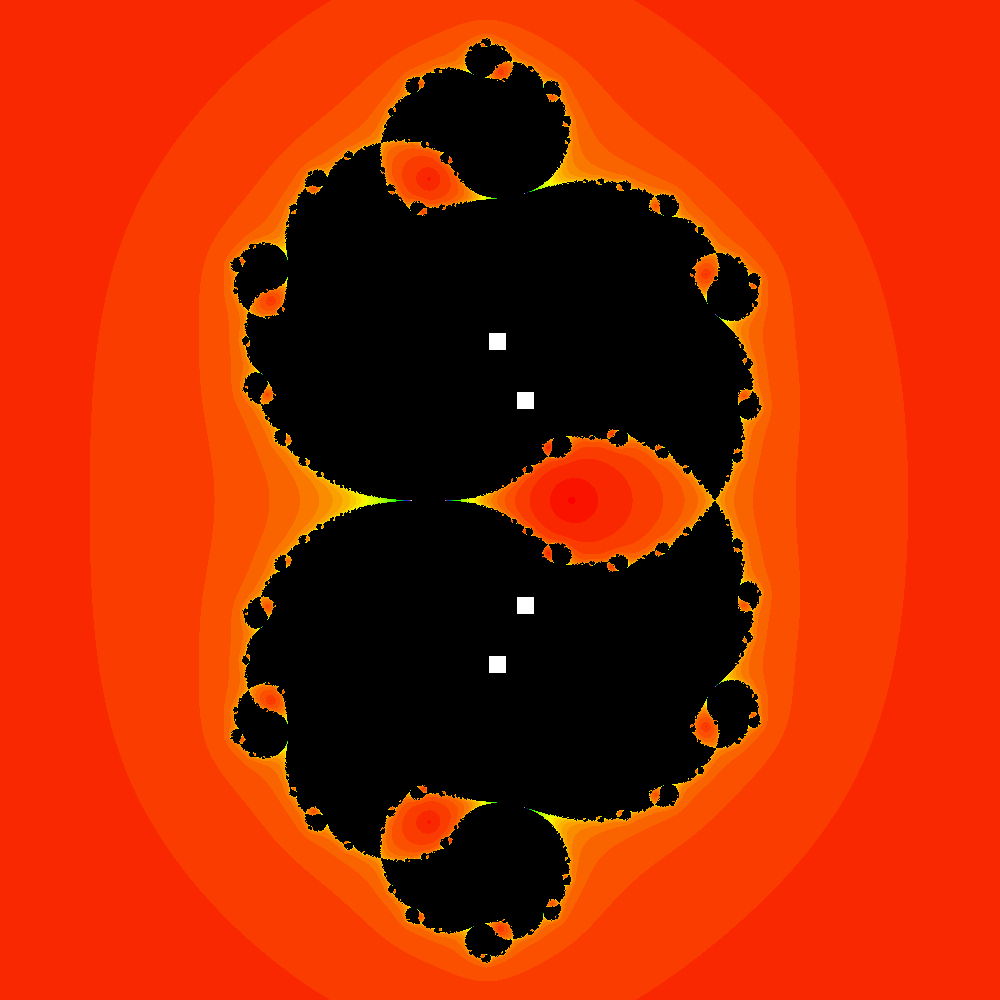};
			\end{axis}
	\end{tikzpicture}}
	\caption{\small{Dynamical planes of the operator \eqref{eq:op4} for different values of $a.$}}
	\label{fig:dinErmakov-O4}
\end{figure}

In Figure \ref{fig:dinErmakov-O4} we show  some dynamical planes for different values of the parameters. We use the same colours as in the dynamical planes of \S~\ref{subsec:op2}. The parameter $a=-20$ lies in the black region which is unbounded to the left for which $z=1$ is attracting. Two free critical points (one up to symmetry) lie in the two petals of the basin of attraction of $z=-1$ (in black) and the other two critical points lie in the basin of attraction of $z=1$ (in green). Similarly,  $a=0.97$ is chosen in the disk described in Figure~\ref{estErmakov} for which $z=1$ is also attracting. The parameter $a=-7$ is chosen on the central red strip of parameters for which the two critical points that do not converge to $-1$ converge to the roots. The parameter $a=9$ lies in the big black component to the right of the parameter plane (see Figure~\ref{fig:paramO4} lower), which corresponds to sets of parameters for which all critical orbits lie in the immediate basin of attraction of $z=-1$.

\subsection{A  sixth order iterative scheme} \label{seccion_tres criticos}

The last example of application of the program we show is a family with three free critical orbits (up to symmetry). In \cite{tresCrit} the authors study the dynamics of a bi-parametric sixth order family of iterative  methods for solving non-linear equations. 

After applying it on quadratic polynomials, for a fixed value of one of the parameters, they obtain the following operator:
{\small
\begin{equation*} \label{op3crit}
 O_a(z)= \frac{z^6 (10-4a+(48-28a+4a^2)z+(69-40a+6a^2)z^2+(56-30a+4a^2)z^3+(28-12a+a^2)z^4+(8-2a)z^5+z^6)}{1+(8-2a)z+(28-12a+a^2)z^2+(56-30a+4a^2)z^3+(69-40a+6a^2)z^4+(48-28a+4a^2)z^5+(10-4a)z^6}, 
\end{equation*}
}
whose derivative is:
{\small
\begin{equation*}\label{dop3crit}
    O^{\prime}_a(z)=\frac{4z^5 (1+z)^4 (-1+(-2+a)z-z^2)P(a,z)}{(1+(8-2a)z+(28-12a+a^2)z^2+(56-30a+4a^2)z^3+(69-40a+6a^2)z^4+(48-28a+4a^2)z^5+(10-4a)z^6)^2}
\end{equation*}
}
where
{\small
\begin{eqnarray*}
P(a,z)&=&-15+6a+(-94+63a-11a^2)z+(-205+170a-49a^2+5a^3)z^2+(-252+206a-56a^2+5a^3)z^3+\\
&+&(-205+170a-49a^2+5a^3)z^4+(-94+63a-11a^2)z^5+(-15+6a)z^6.
\end{eqnarray*}
}
The solutions of equation $-1+(-2+a)z-z^2=0$ are preimages of $z=1$ and the point $z=-1$ is also a preimage of $z=1$. So, there are six free critical points, that are the solutions of the palindromic polynomial of degree six $P(a,z)$. As before, the roots of this polynomial are obtained by doing the change $x=z+\frac{1}{z}$ and solving the polynomial equation of degree three:
\[
-64 +80 a -34 a^2 + 5 a^3 + (-160 + 152 a - 49a^2 + 5 a^3)x + (-94 + 63a- 11a^2) x^2 + (-15 + 6a) x^3 =0.
\]

\begin{figure}[hbt!]
	\centering
	\subfigure{
		\begin{tikzpicture}
			\begin{axis}[width=8cm,  axis equal image, scale only axis,  enlargelimits=false, axis on top]
				\addplot graphics[xmin=3.4,xmax=4.6,ymin=-0.6,ymax=0.6] {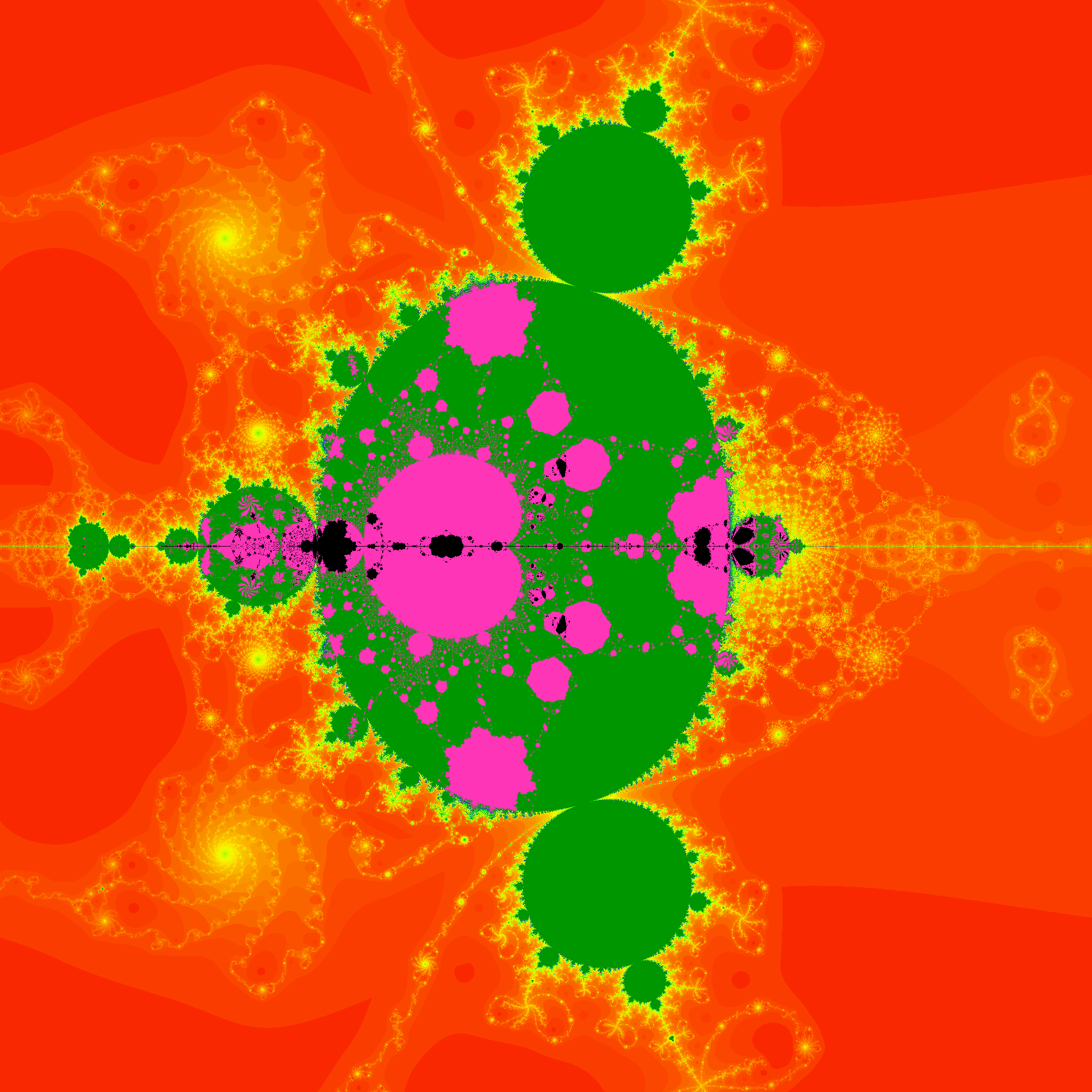};
			\end{axis}
	\end{tikzpicture}}
	\subfigure{
		\begin{tikzpicture}
			\begin{axis}[width=8cm, axis equal image, scale only axis,  enlargelimits=false, axis on top]
				\addplot graphics[xmin=-100,xmax=100,ymin=-100,ymax=100] {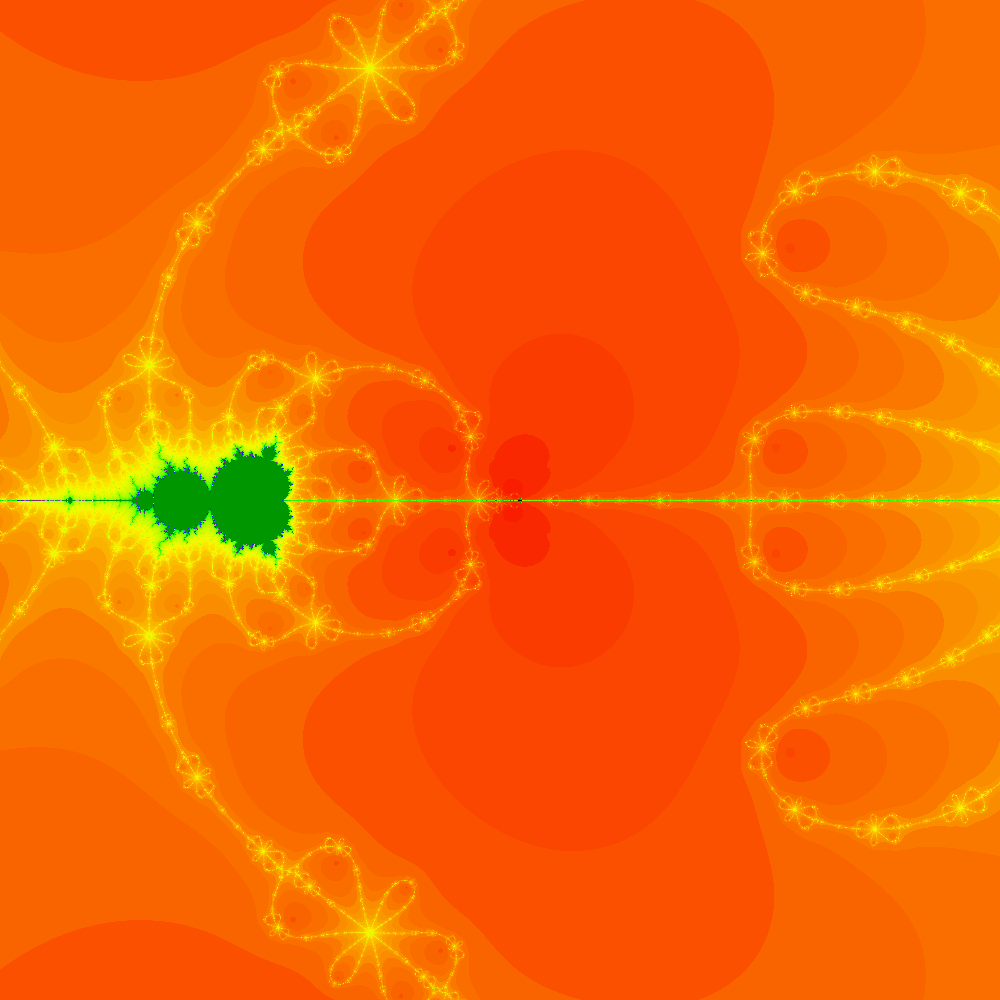};
			\end{axis}
	\end{tikzpicture}}
	\caption{\small{Parameter planes of the operator \eqref{op3crit} obtained using all three free critical orbits simultaneously. Black, pink and green indicate, respectively, that zero, one and two critical orbits have converged to the roots.}}
	\label{fig:param3crit}
\end{figure}

We obtain six different free critical points, only three up to the symmetry of the operator given by $z\rightarrow 1/z$. We plot the parameter plane using simultaneously the three free critical orbits  with independent dynamics (see Figure \ref{fig:param3crit}).  In this case, we use black if no free critical point converges to the roots, pink if one critical orbit converges to a root, green if two critical orbits escape to any of the roots, and a scaling of colours if all critical orbits converge to the roots.

\begin{figure}[h]
\centering
\includegraphics[scale=0.5]{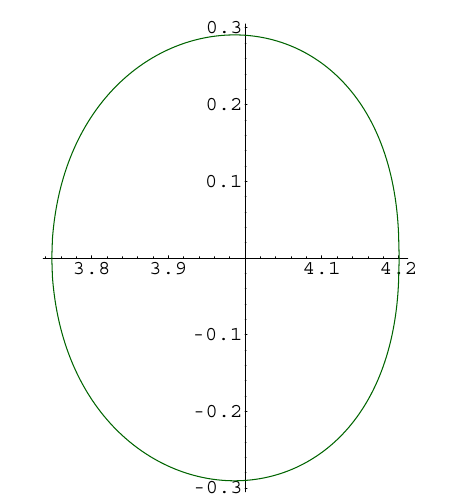}
\caption{Stability region for $z=1$. }\label{est1_trescrit}
\end{figure}

 By checking the values of the parameter where $|O^{\prime}_a(1)|=1$, we have that
the fixed point $z=1$ is attractive inside the curve
\[
1000176-997536\alpha+373160\alpha^2-62056\alpha^3+3871\alpha^4+(124216-62056\alpha+7742\alpha^2)\beta^2+3871\beta^4=0,
\]
being $a=\alpha+i\beta$. This curve can be observed in Figure \ref{est1_trescrit} and corresponds to the biggest green oval (with decorations inside) of Figure \ref{fig:param3crit}).

\begin{figure}[h!]
	\centering
	\subfigure[$a=3.9$]{
		\begin{tikzpicture}
			\begin{axis}[width=6cm,  axis equal image, scale only axis,  enlargelimits=false, axis on top]
				\addplot graphics[xmin=-2,xmax=2,ymin=-2,ymax=2] {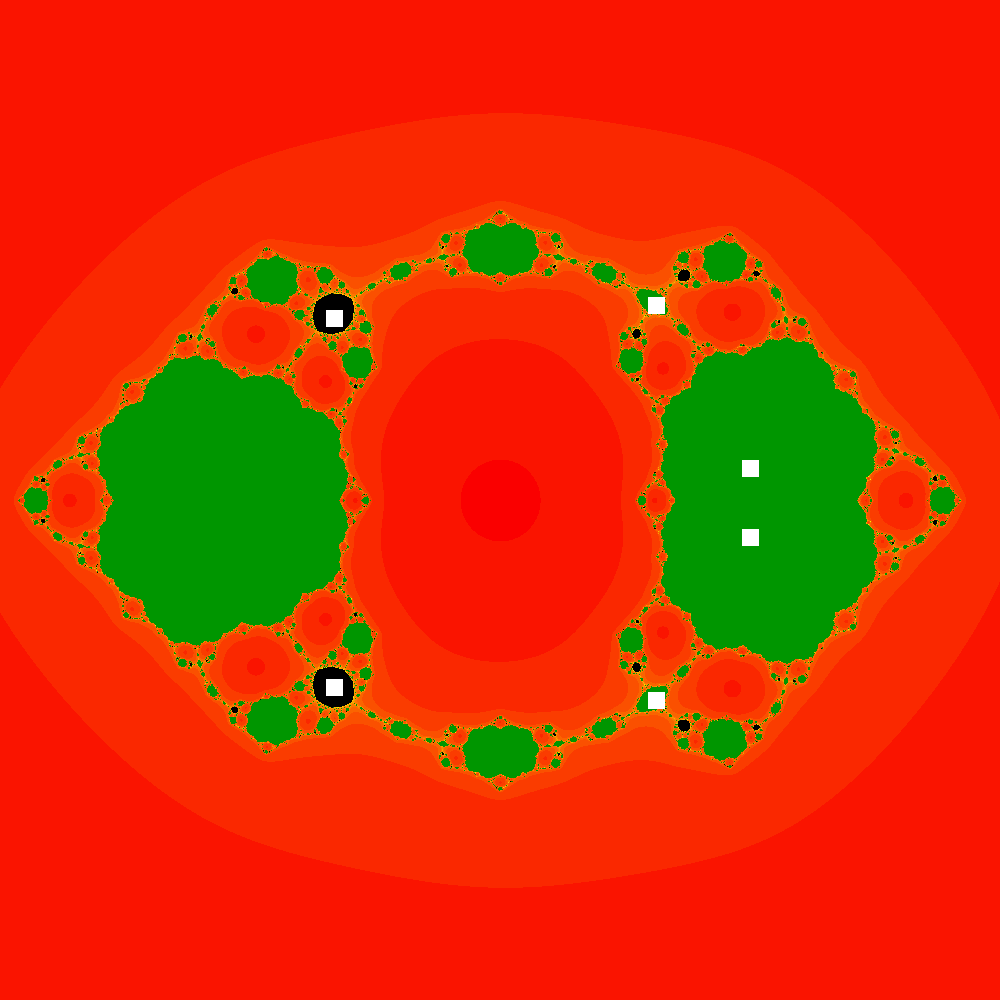};
			\end{axis}
	\end{tikzpicture}}
	\subfigure[$a=3.9+0.04i$]{
		\begin{tikzpicture}
			\begin{axis}[width=6cm,  axis equal image, scale only axis,  enlargelimits=false, axis on top]
				\addplot graphics[xmin=-2,xmax=2,ymin=-2,ymax=2] {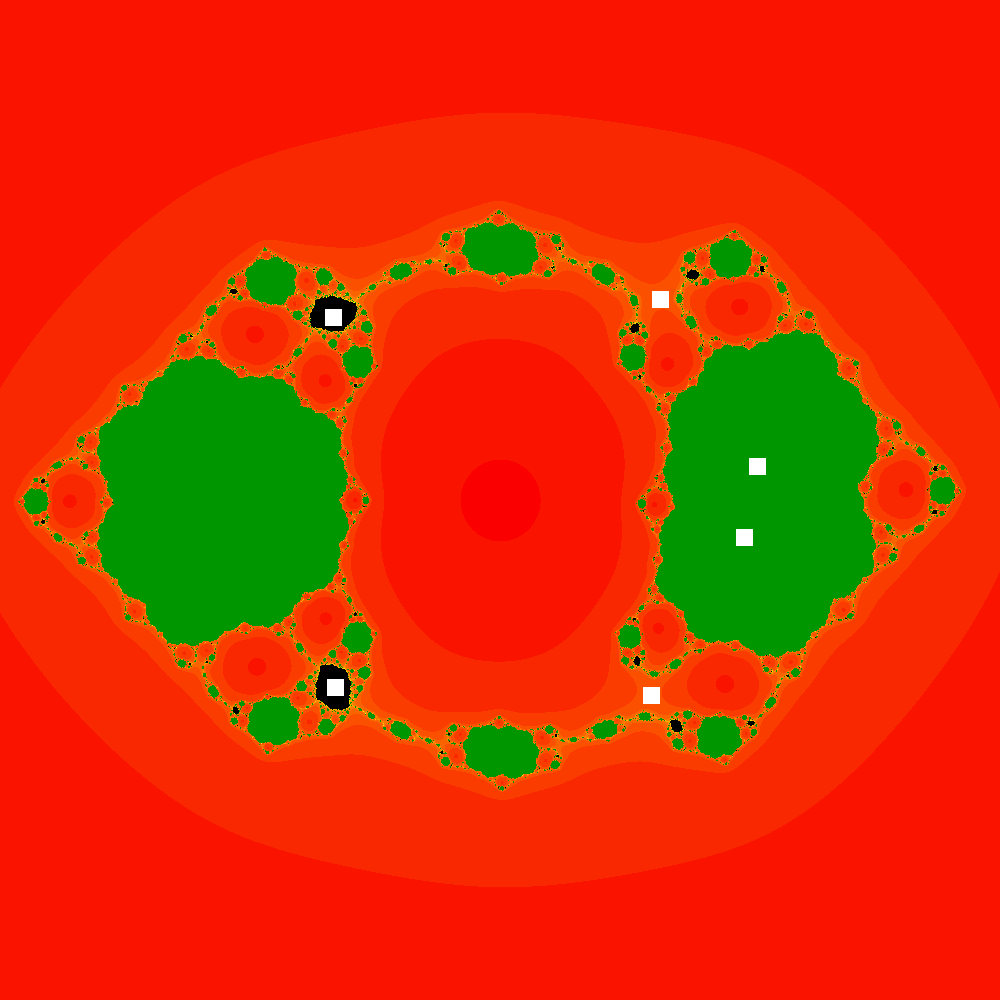};
			\end{axis}
	\end{tikzpicture}}
	
	\subfigure[$a=4.1+0.2i$]{
		\begin{tikzpicture}
			\begin{axis}[width=6cm,  axis equal image, scale only axis,  enlargelimits=false, axis on top]
				\addplot graphics[xmin=-2.25,xmax=2.25,ymin=-2.25,ymax=2.25] {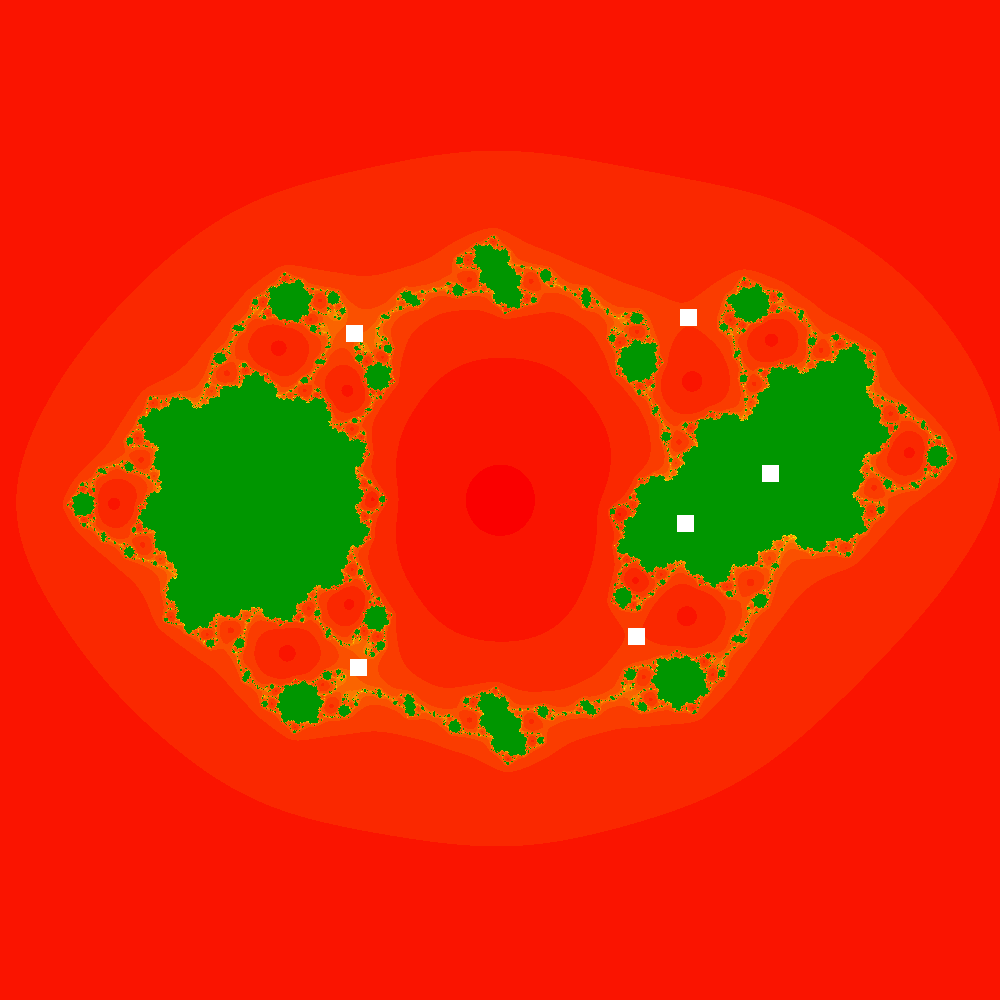};
			\end{axis}
	\end{tikzpicture}}
	\subfigure[$a=4.4+0.4i$]{
		\begin{tikzpicture}
			\begin{axis}[width=6cm,  axis equal image, scale only axis,  enlargelimits=false, axis on top]
				\addplot graphics[xmin=-2.25,xmax=2.75,ymin=-2.5,ymax=2.5] {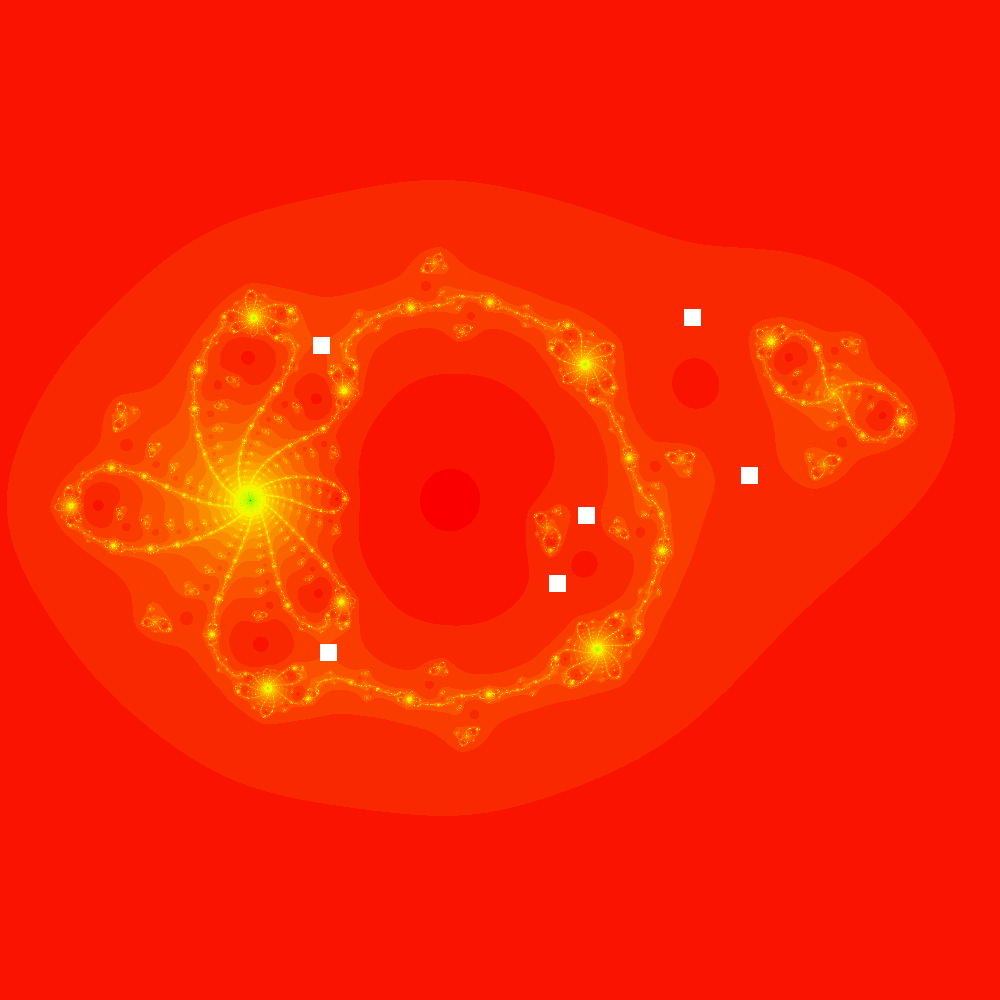};
			\end{axis}
	\end{tikzpicture}}
	\caption{\small{Dynamical planes of the operator \eqref{eq:op4} for different values of $a.$}}
	\label{fig:dynamtrescrit}
\end{figure}

In Figure \ref{fig:dynamtrescrit} we show  some dynamical planes for different values of the parameters. We use the same colours as in the dynamical planes of \S~\ref{subsec:op2}. The first three parameters are chosen within the region for which $z=1$ is attracting. The parameter $a=3.9$ (which appears in black in Figure~\ref{fig:param3crit}) is such that no free critical orbit converges to the roots: four critical points belong to the basin of attraction of $z=1$ (in green) and the other two critical points belong to the basins of attraction of two different attracting fixed points.  The parameter $a=3.9+0.04i$ (which appears in pink in Figure~\ref{fig:param3crit}) is such that two critical orbits belong to the basin of attraction of $z=1$, two critical orbits converge to the roots, and the remaining two critical points converge to the basins of attraction of two different attracting fixed points. The parameter $a=4.1+0.2i$ (which appears in green in Figure~\ref{fig:param3crit}) is such that two critical points belong to the basin of attraction of $z=1$ and the remaining four critical points belong to the basins of attraction of the roots. The last parameter, $a=4.4+0.4i$, is chosen so that all six critical points converge to the roots.

\vspace{0.5cm}
\textbf{Funding}  The first, second and  last authors are supported the project UJI-B2019-18. The first and  last authors are also supported by by the grant PGC2018-095896-B-C22. The second author is also supported by PID2020-118281GB-C32 (MCIU/AEI/FEDER/UE). The third author is supported by Acció 3.2 POSDOC/2020/14 of Universitat Jaume I and by Ayudas Margarita Salas 2021-2023 of Universitat Politècnica de València funded by the Spanish Ministry of Universities (Plan de Recuperación, Transformación y Resiliencia) and European Union-Next generation EU (RD 289/2021 and  UNI/551/2021).

\bibliography{bibliografia}
\bibliographystyle{amsplain}

\end{document}